\documentclass[10pt]{amsart}
\usepackage[T1]{fontenc}
\usepackage{geometry}
\usepackage[latin1] {inputenc}
\usepackage{amsmath}
\usepackage{amsfonts, amssymb, textcomp}
\usepackage[colorlinks=flase, linkcolor=red,urlcolor=green, citecolor=blue]{hyperref}
\usepackage{subeqnarray}

\usepackage{latexsym}
\usepackage{fancyhdr}
\usepackage{longtable}
\usepackage{amsmath, amssymb}
\usepackage{graphicx}
\usepackage{enumerate}

\numberwithin{equation}{section}

\theoremstyle{plain}
\newtheorem{proposition}{Proposition}[section]
\newtheorem{theorem}[proposition]{Theorem}
\newtheorem{lemma}[proposition]{Lemma}
\newtheorem{corollary}[proposition]{Corollary}
\newtheorem{definition}[proposition]{Definition}

\newtheorem{remark}[proposition]{Remark}

\newcommand{\RR}{\mathbb{R}}

\newcommand{\NN}{\mathbb{N}}

\let\on=\operatorname

\newenvironment{demo}[1]{\par\smallskip\noindent{\bf #1.}}{\par\smallskip}

\makeatletter
\@namedef{subjclassname@2020}{%
  \textup{2020} Mathematics Subject Classification}
\makeatother

\newsavebox{\fmbox}
\newenvironment{fmpage}[1]
 {\begin{lrbox}{\fmbox}\begin{minipage}{#1}}
 {\end{minipage}\end{lrbox}\fbox{\usebox{\fmbox}}}

\title[On the equivalence between moderate growth-type conditions]
{On the equivalence between moderate growth-type conditions in the weight matrix setting II}

\author[G.~Schindl]{Gerhard Schindl}

\address{G.~Schindl: Fakult\"at f\"ur Mathematik, Universit\"at Wien, Oskar-Morgenstern-Platz~1, A-1090 Wien, Austria.}
\email{gerhard.schindl@univie.ac.at}

\begin{document}

\begin{abstract}
We continue the study of the known equivalent reformulations of the classical moderate growth condition for weight sequences in the mixed setting; i.e. when dealing with two different sequences. This approach is becoming crucial in the weight matrix setting and also, in particular, when dealing with weight functions in the sense of Braun-Meise-Taylor. It is known that a full generalization to the mixed setting fails, more precisely the condition comparing the growth of the corresponding sequences of quotients and roots is not clear. In the main result we prove a new characterization of this property in terms of the associated weight function; i.e. when the given weight function is based on a weight sequence.
\end{abstract}

\thanks{This research was funded in whole or in part by the Austrian Science Fund (FWF) 10.55776/PAT9445424}
\date{\today}

\maketitle

\section{Introduction}
This article is the direct continuation of \cite{modgrowthstrange} and for the following explanations please see also the introduction there. In the theories of weighted function spaces defined by means of weight sequences $\mathbf{M}=(M_p)_p$, like \emph{ultradifferentiable and ultraholomorphic function spaces} and generalized spaces of \emph{Gelfand-Shilov type,} several basic growth and regularity assumptions on $\mathbf{M}$ appear frequently; see e.g. \cite{Komatsu73}, \cite{BonetMeiseMelikhov07}, \cite{Thilliezdivision}, \cite{optimalflat23}, \cite{nuclearglobal1}, \cite{nuclearglobal2} and the citations therein.\vspace{6pt}

One of the most prominent, classical and important condition is \emph{moderate growth,} in view of \cite{Komatsu73} also known under \emph{(M.2) or stability under ultradifferential operators.} Here, this condition will be abbreviated by \hypertarget{mg}{$(\on{mg})$} and reads as follows:
$$\exists\;C\ge 1\;\forall\;p,q\in\NN:\;M_{p+q}\le C^{p+q+1} M_p M_q.$$
When dealing with weighted spaces it is not restricting the generality to assume $M_0=1$ and in this case one can write $C^{p+q}$ in this condition. \hyperlink{mg}{$(\on{mg})$} implies or even characterizes important and desirable properties of the corresponding weighted function classes.

When $\mathbf{M}$ satisfies mild basic growth conditions then in the literature for \hyperlink{mg}{$(\on{mg})$} several equivalent conditions and reformulations appear and have been used crucially in proofs. For a summary we refer to \cite[Thm. 3.1]{modgrowthstrange} and the citations mentioned there. In particular, it is known that \hyperlink{mg}{$(\on{mg})$} is equivalent to
\begin{equation}\label{mgstrange}
\exists\;A\ge 1\;\forall\;p\in\NN_{>0}:\;\;\;\mu_p\le A(M_p)^{1/p},
\end{equation}
with $\mu_p:=\frac{M_p}{M_{p-1}}$.\vspace{6pt}

The main goal in \cite{modgrowthstrange} has been to study the generalization of this characterizing statement to the mixed setting; i.e. when involving two different sequences. The importance of this question is due to the fact that in \cite{compositionpaper} and \cite{dissertation} we have introduced new (ultradifferentiable) function classes defined in terms of \emph{weight matrices $\mathcal{M}=\{\mathbf{M}^{(x)}: x\in\mathcal{I}\}$,} $\mathcal{I}=(0,+\infty)$ denoting the matrix parameter. We have shown that ultradifferentiable classes defined by a weight function $\omega:[0,+\infty)\rightarrow[0,+\infty)$ in the sense of Braun-Meise-Taylor can be defined by involving the \emph{associated weight matrix} $\mathcal{M}_{\omega}=\{\mathbf{W}^{(x)}: x>0\}$. For the weight function setting we refer to \cite{BraunMeiseTaylor90}, see also \cite{Bjorck66}, and regularity assumptions on $\omega$ are required. Recall that in \cite{BonetMeiseMelikhov07} the weight sequence and weight function approach has been compared; in general both approaches are mutually distinct and for having coincidence \hyperlink{mg}{$(\on{mg})$} is relevant.

Indeed, \hyperlink{mg}{$(\on{mg})$} for some/any $\mathbf{W}^{(x)}$ already implies that both settings are equal and hence this assumption is too restrictive. Similar results are valid for other analogously defined weighted spaces due to the fact that the above comments purely deal with weights and their growth properties.

On the other hand, it is known that the basic assumptions on $\omega$ directly imply a generalization of \hyperlink{mg}{$(\on{mg})$} to the mixed setting, see \eqref{newmoderategrowth}. Note that for weight matrices relevant conditions are naturally arising pair-wise: one considers conditions of \emph{Roumieu- and of Beurling-type,} see e.g. \cite[Sect. 4.1]{compositionpaper}. Inspired by \eqref{newmoderategrowth} we mention the following, see again \cite[Sect. 4.1]{compositionpaper}:

\hypertarget{R-mg}{$(\mathcal{M}_{\{\text{mg}\}})$} \hskip1cm $\forall\;x\in\mathcal{I}\;\exists\;y\in\mathcal{I}\;\exists\;C>0\;\forall\;p,q\in\NN: M^{(x)}_{p+q}\le C^{p+q} M^{(y)}_p M^{(y)}_q$,\par\vskip.3cm

\hypertarget{B-mg}{$(\mathcal{M}_{(\text{mg})})$} \hskip1cm $\forall\;x\in\mathcal{I}\;\exists\;y\in\mathcal{I}\;\exists\;C>0\;\forall\;p,q\in\NN: M^{(y)}_{p+q}\le C^{p+q} M^{(x)}_p M^{(x)}_q$.\par\vskip.3cm

In \cite[Prop. 3.2 \& 3.3, Thm. 3.7, Cor. 3.8]{modgrowthstrange} we have gathered parts of \cite[Thm. 3.1]{modgrowthstrange} which can be transferred to the matrix setting and continued this study in \cite[Prop. 4.2]{modgrowthstrange}; see also the citations in these results. Our attempts to prove a full generalization of \cite[Thm. 3.1]{modgrowthstrange} failed. More precisely, it has turned out that the matrix generalizations of \eqref{mgstrange}, i.e.

\begin{equation}\label{rstrange}
\forall\;x\in\mathcal{I}\;\exists\;y\in\mathcal{I}\;\exists\;A\ge 1\;\forall\;p\in\NN_{>0}:\;\;\;\mu^{(x)}_p\le A(M^{(y)}_p)^{1/p},
\end{equation}
and
\begin{equation}\label{bstrange}
\forall\;x\in\mathcal{I}\;\exists\;y\in\mathcal{I}\;\exists\;A\ge 1\;\forall\;p\in\NN_{>0}:\;\;\;\mu^{(y)}_p\le A(M^{(x)}_p)^{1/p},
\end{equation}
are in general not satisfied for an associated weight matrix. \eqref{rstrange} is called the \emph{quotient-root comparison property of Roumieu-type} and \eqref{bstrange} of \emph{Beurling-type.} Both properties are crucial for technical applications in proofs when working with weight matrices. Indeed, in the main (counter-)example \cite[Thm. 4.8]{modgrowthstrange} we have focused on $\mathcal{M}_{\omega_{\mathbf{M}}}$, where $\omega_{\mathbf{M}}$ denotes the \emph{associated weight function,} see Section \ref{assofctsect}. The relevant information can then be expressed in terms of given $\mathbf{M}$ involving the following new condition \cite[$(4.7)$]{modgrowthstrange}: By \cite[Prop. 4.4]{modgrowthstrange} the matrix $\mathcal{M}_{\omega_{\mathbf{M}}}$ satisfies \eqref{rstrange} and/or \eqref{bstrange} if and only if
\begin{equation}\label{genmgintro}
\exists\;d\in\NN_{>0}\;\exists\;A\ge 1\;\forall\;p\in\NN_{>0}:\;\;\;\mu_p\le A(M_{dp})^{\frac{1}{dp}}.
\end{equation}
When $p\mapsto(M_p)^{1/p}$ is nondecreasing, which is a standard assumption in this context, then \eqref{genmgintro} is weaker than \eqref{mgstrange} and \eqref{genmgintro} also motivates to consider the new growth index $g(\mathbf{M})$ which is defined to be the minimal $d$ such that the estimate in \eqref{genmgintro} is valid. However, in \cite{modgrowthstrange} we have not been able to prove a characterization of \eqref{genmgintro} in terms of $\omega_{\mathbf{M}}$ and this is the main content of this work. Such a characterization is natural and desirable when working in the weight function setting since it directly gives a growth restriction on the (associated) weight function. In order to proceed we are applying new information concerning the \emph{generalized lower Legendre conjugate (or envelope)} for weight functions obtained in the recent work \cite{genLegendreconj}; see Theorem \ref{multthm}. This gives a direct application of this conjugate and underlines its importance in a concrete direction. We formulate the main results in a general mixed setting and give also a connection to the notion of \emph{weak separativity} considered in \cite{matsumoto}. Therefore, also the technical preparatory results admit applications in different contexts and independently we can strengthen some results from \cite{modgrowthstrange}; e.g. obtain as a by-product of these statements that \eqref{genmgintro} is preserved under equivalence of weight sequences. Note also that the results for $\omega_{\mathbf{M}}$ can be transferred to abstractly given weights $\omega$ in the sense of Braun-Meise-Taylor when replacing $\mathbf{M}$ by $\mathbf{W}^{(1)}$; see \cite[Prop. 4.7]{modgrowthstrange}.\vspace{6pt}

We expect that the quotient-root comparison properties and hence the characterizations shown in this work (and in \cite{modgrowthstrange}) admit direct applications to different weighted settings. Indeed, the results \cite[Cor. 9 \& 10, Thm. 11]{whitneyextensionmixedweightfunctionII} have been inspiring for this research, see also comments $(i)-(iv)$ just after \cite[Thm. 4.8]{modgrowthstrange}. Very recently, condition \eqref{mgstrange} has been applied in the crucial technical result \cite[Lemma 2]{CordaroFuerdoes24} and thus using the quotient-root comparison properties \eqref{rstrange}, \eqref{bstrange} together with the characterizations established in this work (and in \cite{modgrowthstrange}) gives the idea to transfer the main statements shown in \cite{CordaroFuerdoes24} from the weight sequence to the weight function setting.\vspace{6pt}

The paper is structured as follows: In Section \ref{weightscondsect} we gather all relevant information concerning weight sequences and (associated) weight functions. In Section \ref{auxsection} we revisit some technical auxiliary sequences needed in the proofs of the main results. We also prove new estimates for the moderate growth index $g(\cdot)$ for sequences belonging to $\mathcal{M}_{\omega}$; see Lemma \ref{technlemma1}. In Section \ref{preparatorysection} we obtain new characterizations of \eqref{rstrange} resp. \eqref{bstrange} in terms of these auxiliary sequences, see Propositions \ref{technlemma2}, \ref{technlemma3} and \ref{technlemma3consequ}. Using this we can show that \eqref{genmgintro} is preserved under equivalence of weight sequences; see Theorem \ref{technlemma4} and Remark \ref{technlemma4rem}. In Section \ref{mainsection} we prove the main statements of this article: Theorem \ref{upptrafothm} and Corollaries \ref{omegaMcharacterzing} and \ref{omegaMcharacterzing1}. Corollary \ref{omegaMcharacterzing1} gives the desired characterization of \eqref{genmgintro} in terms of $\omega_{\mathbf{M}}$. In Sections \ref{preparatorysectionweightfct} and \ref{mainsectweightfct} comments on abstractly given weight functions in the sense of Braun-Meise-Taylor are given. Finally, in Section \ref{counterexsection} we are concerned with the problem if \eqref{rstrange} resp. \eqref{bstrange} can be assumed w.l.o.g. when switching to an equivalent weight matrix. However, a general answer to this question seems to be difficult since we can only derive a necessary condition but which is too general; see Proposition \ref{mainmgthm3} and Lemma \ref{mainmgthm3problemlemma}.\vspace{6pt}

\textbf{Acknowledgements.} The author of this article thanks the two anonymous referees for the careful reading and the valuable suggestions.

\section{Weights and conditions}\label{weightscondsect}
\subsection{General notation}
We write $\NN:=\{0,1,2,\dots\}$, $\NN_{>0}:=\{1,2,3,\dots\}$ and $\RR_{>0}:=(0,+\infty)$.

\subsection{Weight sequences}\label{weightsequences}
Given a sequence $\mathbf{M}=(M_p)_p\in\RR_{>0}^{\NN}$ we also use the notation $\mu=(\mu_p)_p$ with $\mu_p:=\frac{M_p}{M_{p-1}}$, $p\ge 1$, $\mu_0:=1$, and analogously for all other arising sequences. $\mathbf{M}$ is called \emph{normalized,} if $1=M_0\le M_1$ is valid. For the point-wise product we simply write $\mathbf{M}\mathbf{N}$ and for the $\ell$-th power of $\mathbf{M}$ we write $\mathbf{M}^{\ell}$, $\ell>0$. This notation should not be mixed with sequences having index $\ell>0$ and belonging to a certain given weight matrix; here we write $\mathbf{M}^{(\ell)}$.\vspace{6pt}

$\mathbf{M}$ is called \emph{log-convex} if
$$\forall\;p\in\NN_{>0}:\;M_p^2\le M_{p-1} M_{p+1},$$
equivalently if $p\mapsto\mu_p$ is nondecreasing. If $\mathbf{M}$ is log-convex and normalized, then $\mu_p\ge 1$ for all $p\in\NN$, $p\mapsto M_p$ and $p\mapsto(M_p)^{1/p}$ are nondecreasing, $(M_p)^{1/p}\le\mu_p$ for all $p\in\NN_{>0}$ and finally $M_pM_q\le M_{p+q}$ for all $p,q\in\NN$; see e.g. \cite[Lemmas 2.0.4 \& 2.0.6]{diploma}.\vspace{6pt}

For our purpose it is convenient to consider the following set of sequences
$$\hypertarget{LCset}{\mathcal{LC}}:=\{\mathbf{M}\in\RR_{>0}^{\NN}:\;\mathbf{M}\;\text{is normalized, log-convex},\;\lim_{p\rightarrow+\infty}(M_p)^{1/p}=+\infty\}.$$
We see that $\mathbf{M}\in\hyperlink{LCset}{\mathcal{LC}}$ if and only if $1=\mu_0\le\mu_1\le\dots$ and $\lim_{p\rightarrow+\infty}\mu_p=+\infty$ (see e.g. \cite[p. 104]{compositionpaper}) and there is a one-to-one correspondence between $\mathbf{M}$ and $\mu=(\mu_p)_p$ by taking $M_p:=\prod_{i=0}^p\mu_i$.\vspace{6pt}

Let $\mathbf{M},\mathbf{N}\in\RR_{>0}^{\NN}$ be given, we write $\mathbf{M}\le\mathbf{N}$ if $M_p\le N_p$ for all $p\in\NN$ and $\mathbf{M}\hypertarget{preceq}{\preceq}\mathbf{N}$ if $\sup_{p\in\NN_{>0}}\left(\frac{M_p}{N_p}\right)^{1/p}<+\infty.$ We call $\mathbf{M}$ and $\mathbf{N}$ \emph{equivalent,} formally denoted by $\mathbf{M}\hypertarget{approx}{\approx}\mathbf{N}$, if $\mathbf{M}\hyperlink{preceq}{\preceq}\mathbf{N}$ and $\mathbf{N}\hyperlink{preceq}{\preceq}\mathbf{M}$. Property \hyperlink{mg}{$(\on{mg})$} is clearly preserved under \hyperlink{approx}{$\approx$}.

\subsection{Associated weight function}\label{assofctsect}
Let $\mathbf{M}\in\RR_{>0}^{\NN}$ (with $M_0=1$), then the \emph{associated function} $\omega_{\mathbf{M}}: \RR_{\ge 0}\rightarrow\RR\cup\{+\infty\}$ is defined by
\begin{equation*}\label{assofunc}
\omega_{\mathbf{M}}(t):=\sup_{p\in\NN}\log\left(\frac{t^p}{M_p}\right)\;\;\;\text{for}\;t\ge 0,
\end{equation*}
with the conventions $0^0:=1$ and $\log(0)=-\infty$. This ensures $\omega_{\mathbf{M}}(0)=0$ and $\omega_{\mathbf{M}}(t)\ge 0$ for any $t\ge 0$ since $\frac{t^0M_0}{M_0}=1$ for all $t\ge 0$. For an abstract introduction of the associated function we refer to the classical work \cite[Chapitre I]{mandelbrojtbook}, see also \cite[Definition 3.1]{Komatsu73} and to the more recent article \cite{regseqpaper}. It is immediate by definition that $\mathbf{M}\le\mathbf{N}$ implies $\omega_{\mathbf{N}}(t)\le\omega_{\mathbf{M}}(t)$ for all $t\ge 0$ and we have the following identity for the $\ell$-th power $\mathbf{M}^{\ell}$ of $\mathbf{M}$:
\begin{equation}\label{powersub}
\forall\;\ell>0\;\forall\;t\ge 0:\;\;\;\omega_{\mathbf{M}^{\ell}}(t)=\ell\omega_{\mathbf{M}}(t^{1/\ell}).
\end{equation}

We define the \emph{counting function} $\Sigma_{\mathbf{M}}$ by
\begin{equation}\label{counting}
\Sigma_{\mathbf{M}}(t):=|\{p\in\NN_{>0}:\;\;\;\mu_p\le t\}|,\;\;\;t\ge 0.
\end{equation}
It is known that for given $\mathbf{M}\in\hyperlink{LCset}{\mathcal{LC}}$ the functions $\omega_{\mathbf{M}}$ and $\Sigma_{\mathbf{M}}$ are related by the following integral representation formula; see \cite[1.8. III]{mandelbrojtbook}, \cite[$(3.11)$]{Komatsu73} and also \cite[Lemma 2.4]{regseqpaper}:
\begin{equation}\label{intrepr}
\forall\;t\ge 0:\;\;\;\omega_{\mathbf{M}}(t)=\int_0^t\frac{\Sigma_{\mathbf{M}}(u)}{u}du=\int_{\mu_1}^t\frac{\Sigma_{\mathbf{M}}(u)}{u}du.
\end{equation}
Consequently, $\omega_{\mathbf{M}}$ vanishes on $[0,\mu_1]$, in particular on the unit interval.

Finally, if $\mathbf{M}\in\hyperlink{LCset}{\mathcal{LC}}$, or even when $\mathbf{M}$ is log-convex, $\lim_{p\rightarrow+\infty}(M_p)^{1/p}=+\infty$ and $M_0=1$, then we can compute $\mathbf{M}$ by involving $\omega_{\mathbf{M}}$ as follows; see \cite[Chapitre I, 1.4, 1.8]{mandelbrojtbook}, \cite[Prop. 3.2]{Komatsu73} and also \cite[Sect. 2 \& 3]{regseqpaper}:
\begin{equation}\label{Prop32Komatsu}
M_p=\sup_{t\ge 0}\frac{t^p}{\exp(\omega_{\mathbf{M}}(t))},\;\;\;p\in\NN.
\end{equation}

\subsection{Weight matrices}
Let $\mathcal{I}=\RR_{>0}$ be the index set. A \emph{weight matrix} $\mathcal{M}$ associated with $\mathcal{I}$ is the set $\mathcal{M}=\{\mathbf{M}^{(x)}: x>0\}$ such that $\mathbf{M}^{(x)}\le\mathbf{M}^{(y)}$ for all $0<x\le y$, where $\mathbf{M}^{(x)}$ is a weight sequence for all $x\in\mathcal{I}$. $\mathcal{M}$ is called
\begin{itemize}
\item[$(*)$]\emph{constant} if $\mathbf{M}^{(x)}\hyperlink{approx}{\approx}\mathbf{M}^{(y)}$ for all $x,y>0$ and

\item[$(*)$] \emph{standard log-convex} if $\mathbf{M}^{(x)}\in\hyperlink{LCset}{\mathcal{LC}}$ for all $x>0$.
\end{itemize}

Let $\mathcal{M}=\{\mathbf{M}^{(x)}: x>0\}$ and $\mathcal{N}=\{\mathbf{N}^{(y)}: y>0\}$ be given. We write $\mathcal{M}\hypertarget{Mroumpreceq}{\{\preceq\}}\mathcal{N}$ if
$$\forall\;x>0\;\exists\;y>0:\;\;\;\mathbf{M}^{(x)}\hyperlink{preceq}{\preceq}\mathbf{N}^{(y)},$$
and $\mathcal{M}\hypertarget{Mbeurpreceq}{(\preceq)}\mathcal{N}$ if
$$\forall\;x>0\;\exists\;y>0:\;\;\;\mathbf{M}^{(y)}\hyperlink{preceq}{\preceq}\mathbf{N}^{(x)}.$$
$\mathcal{M}$ and $\mathcal{N}$ are \emph{R-equivalent,} written $\mathcal{M}\{\approx\}\mathcal{N}$, if
$\mathcal{M}\hyperlink{Mroumpreceq}{\{\preceq\}}\mathcal{N}$ and $\mathcal{N}\hyperlink{Mroumpreceq}{\{\preceq\}}\mathcal{M}$ and \emph{B-equivalent,} denoted by $\mathcal{M}(\approx)\mathcal{N}$, if $\mathcal{M}\hyperlink{Mbeurpreceq}{(\preceq)}\mathcal{N}$ and $\mathcal{N}\hyperlink{Mbeurpreceq}{(\preceq)}\mathcal{M}$.

\subsection{Weight functions and associated weight matrices}\label{assomatrixsection}
$\omega:[0,+\infty)\rightarrow[0,+\infty)$ is called a \emph{weight function} (in the terminology of \cite[Section 2.1]{index} and \cite[Section 2.2]{sectorialextensions}), if it is continuous, nondecreasing, $\omega(0)=0$ and $\lim_{t\rightarrow+\infty}\omega(t)=+\infty$. If $\omega$ satisfies in addition $\omega(t)=0$ for all $t\in[0,1]$, then $\omega$ is called \emph{normalized.} For convenience we write that $\omega$ has $\hypertarget{om0}{(\omega_0)}$ if it is a normalized weight.\vspace{6pt}

Moreover we consider the following conditions; this list of properties has already been used in ~\cite{dissertation}.

\begin{itemize}
\item[\hypertarget{om1}{$(\omega_1)}$] $\omega(2t)=O(\omega(t))$ as $t\rightarrow+\infty$, i.e. $\exists\;L\ge 1\;\forall\;t\ge 0:\;\;\;\omega(2t)\le L(\omega(t)+1)$.


\item[\hypertarget{om3}{$(\omega_3)$}] $\log(t)=o(\omega(t))$ as $t\rightarrow+\infty$.

\item[\hypertarget{om4}{$(\omega_4)$}] $\varphi_{\omega}:t\mapsto\omega(e^t)$ is a convex function on $\RR$.


\item[\hypertarget{om6}{$(\omega_6)$}] $\exists\;H\ge 1\;\forall\;t\ge 0:\;2\omega(t)\le\omega(H t)+H$.
\end{itemize}

Finally, we recall the \emph{strong nonquasianalyticity} condition for weight functions
\begin{equation}\label{assostrongnq}
\exists\;C\ge 1\;\forall\;y\ge 0:\;\;\;\int_1^{+\infty}\frac{\omega(yt)}{t^2}dt\le C\omega(y)+C.
\end{equation}

For convenience we define the set
$$\hypertarget{omset0}{\mathcal{W}_0}:=\{\omega:[0,\infty)\rightarrow[0,\infty): \omega\;\text{has}\;\hyperlink{om0}{(\omega_0)},\hyperlink{om3}{(\omega_3)},\hyperlink{om4}{(\omega_4)}\}.$$
For any $\omega\in\hyperlink{omset0}{\mathcal{W}_0}$ we define the \emph{Legendre-Fenchel-Young-conjugate} of $\varphi_{\omega}$ by
\begin{equation}\label{legendreconjugate}
\varphi^{*}_{\omega}(x):=\sup\{x y-\varphi_{\omega}(y): y\ge 0\},\;\;\;x\ge 0,
\end{equation}
with the following properties, see e.g. \cite[Rem. 1.3, Lemma 1.5]{BraunMeiseTaylor90}: It is convex and nondecreasing, $\varphi^{*}_{\omega}(0)=0$, $\varphi^{**}_{\omega}=\varphi_{\omega}$, $\lim_{x\rightarrow+\infty}\frac{x}{\varphi^{*}_{\omega}(x)}=0$ and finally $x\mapsto\frac{\varphi_{\omega}(x)}{x}$ and $x\mapsto\frac{\varphi^{*}_{\omega}(x)}{x}$ are nondecreasing on $[0,+\infty)$. Note that by normalization we can extend the supremum in \eqref{legendreconjugate} from $y\ge 0$ to $y\in\RR$ without changing the value of $\varphi^{*}_{\omega}(x)$ for given $x\ge 0$.\vspace{6pt}

We recall the following known result, see e.g. \cite[Lemma 2.8]{testfunctioncharacterization} resp. \cite[Lemma 2.4]{sectorialextensions} and the references mentioned in the proofs there.

\begin{lemma}\label{assoweightomega0}
Let $\mathbf{M}\in\hyperlink{LCset}{\mathcal{LC}}$, then $\omega_{\mathbf{M}}\in\hyperlink{omset0}{\mathcal{W}_0}$ holds and \hyperlink{om6}{$(\omega_6)$} for $\omega_{\mathbf{M}}$ if and only if $\mathbf{M}$ has \hyperlink{mg}{$(\on{mg})$}.
\end{lemma}

Let $\sigma,\tau$ be weight functions, we write $\sigma\hypertarget{ompreceq}{\preceq}\tau$ if $\tau(t)=O(\sigma(t))\;\text{as}\;t\rightarrow+\infty$
and call the weights equivalent, denoted by $\sigma\hypertarget{sim}{\sim}\tau$, if
$\sigma\hyperlink{ompreceq}{\preceq}\tau$ and $\tau\hyperlink{ompreceq}{\preceq}\sigma$.

All listed properties, except the convexity condition \hyperlink{om4}{$(\omega_4)$}, are clearly preserved under \hyperlink{sim}{$\sim$}.

\begin{remark}\label{omega6rem}
\emph{Let $\omega:[0,+\infty)\rightarrow[0,+\infty)$ be nondecreasing and $\lim_{t\rightarrow+\infty}\omega(t)=+\infty$; i.e. $\omega$ is a weight function in the sense of \cite{genLegendreconj} (and of \cite{index}). Then the following are equivalent:}
\begin{itemize}
\item[$(a)$] \emph{$\omega$ satisfies \hyperlink{om6}{$(\omega_6)$}.}

\item[$(b)$] \emph{$\omega$ satisfies}
$$\exists\;a>1\;\exists\;H\ge 1\;\forall\;t\ge 0:\;\;\;a\omega(t)\le\omega(Ht)+H.$$

\item[$(c)$] \emph{$\omega$ satisfies}
$$\forall\;a>1\;\exists\;H\ge 1\;\forall\;t\ge 0:\;\;\;a\omega(t)\le\omega(Ht)+H.$$
\end{itemize}
\emph{$(a)\Rightarrow(b)$ and $(c)\Rightarrow(a)$ are clear (take $a:=2$). For $(b)\Rightarrow(c)$ let $b>1$ be arbitrary, then choose $n\in\NN_{>0}$ such that $a^n\ge b$ with $a$ appearing in $(b)$, and iterate the estimate $n$-times.}
\end{remark}

We summarize some facts which are shown in \cite[Sect. 5]{compositionpaper} and \cite[Sect. 4 \& 5]{dissertation} and are needed in this work. All properties listed below are valid for $\omega\in\hyperlink{omset0}{\mathcal{W}_0}$.

\begin{itemize}
\item[$(i)$] With each $\omega\in\hyperlink{omset0}{\mathcal{W}_0}$ we can associate a standard log-convex weight matrix $\mathcal{M}_{\omega}:=\{\mathbf{W}^{(\ell)}=(W^{(\ell)}_p)_{p\in\NN}: \ell>0\}$ by\vspace{6pt}

    \centerline{$W^{(\ell)}_p:=\exp\left(\frac{1}{\ell}\varphi^{*}_{\omega}(\ell p)\right)$.}\vspace{6pt}

For the corresponding sequence of quotients write $\vartheta^{(\ell)}=(\vartheta^{(\ell)}_p)_{p\in\NN}$, i.e. $\vartheta^{(\ell)}_p:=\frac{W^{(\ell)}_p}{W^{(\ell)}_{p-1}}$ for $p\in\NN_{>0}$ and $\vartheta^{(\ell)}_0:=1$. Note that the sequences are even ordered w.r.t. their quotients which means $\vartheta^{(\ell_1)}\le\vartheta^{(\ell_2)}$ for any $0<\ell_1\le\ell_2$; see \cite[Sect. 2.5]{whitneyextensionweightmatrix}.

\item[$(ii)$] $\mathcal{M}_{\omega}$ satisfies
    \begin{equation}\label{newmoderategrowth}
    \forall\;\ell>0\;\forall\;p,q\in\NN:\;\;\;W^{(\ell)}_{p+q}\le W^{(2\ell)}_pW^{(2\ell)}_q,
    \end{equation}
    which implies both \hyperlink{R-mg}{$(\mathcal{M}_{\{\on{mg}\}})$} and \hyperlink{B-mg}{$(\mathcal{M}_{(\on{mg})})$}.

\item[$(iii)$] \hyperlink{om6}{$(\omega_6)$} holds if and only if some/each $\mathbf{W}^{(\ell)}$ satisfies \hyperlink{mg}{$(\on{mg})$} if and only if $\mathbf{W}^{(\ell_1)}\hyperlink{approx}{\approx}\mathbf{W}^{(\ell_2)}$ for each $\ell_1,\ell_2>0$. Consequently \hyperlink{om6}{$(\omega_6)$} is characterizing the situation when $\mathcal{M}_{\omega}$ is constant.

\item[$(iv)$] We have $\omega\hyperlink{sim}{\sim}\omega_{\mathbf{W}^{(\ell)}}$ for each $\ell>0$, more precisely
\begin{equation}\label{goodequivalenceclassic}
\forall\;\ell>0\,\,\exists\,D_{\ell}>0\;\forall\;t\ge 0:\;\;\;\ell\omega_{\mathbf{W}^{(\ell)}}(t)\le\omega(t)\le 2\ell\omega_{\mathbf{W}^{(\ell)}}(t)+D_{\ell};
\end{equation}
see \cite[Theorem 4.0.3, Lemma 5.1.3]{dissertation} and also \cite[Lemma 2.5]{sectorialextensions}. Note that for having \eqref{goodequivalenceclassic} the convexity condition \hyperlink{om4}{$(\omega_4)$} is indispensable; more precisely this property is required for the second estimate in \eqref{goodequivalenceclassic}.

\item[$(v)$] The following formulas are immediate by definition:

\begin{equation}\label{transform}
\forall\;x>0\;\forall\;\ell\in\NN_{>0}\;\forall\;p\in\NN:\;\;\;W^{(\ell x)}_p=(W^{(x)}_{\ell p})^{\frac{1}{\ell}},
\end{equation}
and
\begin{equation}\label{transform1}
\forall\;x>0\;\forall\;\ell\in\NN_{>0}\;\forall\;p\in\NN:\;\;\;W^{(x/\ell)}_{p\ell}=(W^{(x)}_p)^{\ell}.
\end{equation}
\end{itemize}

Let $\mathbf{M}\in\hyperlink{LCset}{\mathcal{LC}}$ be given, then $\omega_{\mathbf{M}}\in\hyperlink{omset0}{\mathcal{W}_0}$, see Lemma \ref{assoweightomega0}, and so it makes sense to define and study the matrix
$$\mathcal{M}_{\omega_{\mathbf{M}}}:=\{\mathbf{M}^{(\ell)}: \ell>0\}.$$
In this case, let us recall \cite[$(2.11)$, p. 9]{modgrowthstrange}:
\begin{equation}\label{transformM}
\forall\;p\in\NN:\;\;\;M_p=M^{(1)}_p.
\end{equation}

\subsection{New characterizing condition}
The following new crucial condition for $\mathbf{M}\in\hyperlink{LCset}{\mathcal{LC}}$, or even for $\mathbf{M}\in\RR_{>0}^{\NN}$, has been introduced in \cite{modgrowthstrange}; see \cite[$(4.7)$]{modgrowthstrange}: For fixed $d\in\NN_{>0}$ consider
\begin{equation}\label{genmg}
\exists\;A\ge 1\;\forall\;p\in\NN_{>0}:\;\;\;\mu_p\le A(M_{dp})^{\frac{1}{dp}}.
\end{equation}
Based on this condition, in \cite[p. 24]{modgrowthstrange} we have given the following definition:

\begin{definition}\label{gindexdef}
\emph{Let $\mathbf{M}\in\hyperlink{LCset}{\mathcal{LC}}$, then the \emph{moderate growth index} $g(\mathbf{M})$ is defined by}
$$g(\mathbf{M}):=\min\{d\in\NN_{>0}:\;\;\;\eqref{genmg}\;\emph{is valid}\},$$
\emph{and let us set $g(\mathbf{M}):=+\infty$ if \eqref{genmg} is violated.}
\end{definition}

So for $\mathbf{M}\in\hyperlink{LCset}{\mathcal{LC}}$ we have that $g(\mathbf{M})=1$ holds if and only if $\mathbf{M}$ has \hyperlink{mg}{$(\on{mg})$} and $g(\mathbf{M})<+\infty$ if and only if $\mathcal{M}_{\omega_{\mathbf{M}}}$ has \eqref{rstrange} and/or \eqref{bstrange} by the main statement \cite[Prop. 4.4]{modgrowthstrange}. And, by \cite[Prop. 4.7]{modgrowthstrange}, the analogous result holds for $\mathcal{M}_{\omega}$ with abstractly given $\omega\in\hyperlink{omset0}{\mathcal{W}_0}$ when $\mathbf{M}$ is replaced by $\mathbf{W}^{(1)}$.

\section{Technical auxiliary sequences}\label{auxsection}
Before we start proving the main results, first we have to recall some useful notation. For any $\mathbf{M}\in\RR_{>0}^{\NN}$ and $a\in\NN_{>0}$ let us introduce the sequence $\widetilde{\mathbf{M}}^a$ given by
\begin{equation}\label{tildeweightsequences}
\widetilde{\mathbf{M}}^a_p:=(M_{ap})^{1/a}.
\end{equation}
This notation has already been used with $a:=4$ in \cite[$(3.2)$, Sect. 3]{modgrowthstrange} and for these auxiliary sequences we refer also to \cite[Lemma 2.2]{subaddlike}, \cite[Lemma 6.5]{PTTvsmatrix} and \cite[Sect. 3.3, $(3.26)$]{weightedentireinclusion1}. We denote by $\widetilde{\mu}^a=(\widetilde{\mu}^a_p)_{p\in\NN}$ the corresponding sequence of quotients. Set again $\widetilde{\mu}^a_0:=1$ and by definition
\begin{equation}\label{tildeweightsequences1}
\widetilde{\mu}^a_p=\left(\frac{M_{ap}}{M_{a(p-1)}}\right)^{1/a}=(\mu_{ap-a+1}\cdots\mu_{ap})^{1/a},\;\;\;p\in\NN_{>0}.
\end{equation}
Consequently, one has the following properties:
\begin{itemize}
\item[$(*)$] $\widetilde{(\widetilde{\mathbf{M}}^a)}^b=\widetilde{\mathbf{M}}^{ab}$ for any $a,b\in\NN_{>0}$.

\item[$(*)$] $\lim_{p\rightarrow+\infty}(M_p)^{1/p}=+\infty$ if and only if $\lim_{p\rightarrow+\infty}(\widetilde{M}^a_p)^{1/p}=+\infty$ for some/any $a\in\NN_{>0}$.

\item[$(*)$] If $M_0=1$ and $\mathbf{M}$ is log-convex, then $\widetilde{\mathbf{M}}^a\le\widetilde{\mathbf{M}}^b$ for all $a\le b$ and

\item[$(*)$] by \eqref{tildeweightsequences1} we get (see also \cite[$(4.9)$]{modgrowthstrange}):
\begin{equation}\label{tildeweightsequences2}
\forall\;a\in\NN_{>0}\;\forall\;p\in\NN_{>0}:\;\;\;\mu_{a(p-1)}\le\mu_{ap-a+1}\le\widetilde{\mu}^a_p\le\mu_{ap}.
\end{equation}

\item[$(*)$] If $\mathbf{M}\in\hyperlink{LCset}{\mathcal{LC}}$, then $\widetilde{\mathbf{M}}^a\hyperlink{preceq}{\preceq}\mathbf{M}$ for some/any $a\in\NN_{>0}$ if and only if $\mathbf{M}$ satisfies \hyperlink{mg}{$(\on{mg})$} (see the proof of \cite[Lemma 2.2]{subaddlike}).

\item[$(*)$] If $\mathbf{M}\in\hyperlink{LCset}{\mathcal{LC}}$, then $\widetilde{\mathbf{M}}^a\in\hyperlink{LCset}{\mathcal{LC}}$ for any $a\in\NN_{>0}$.
\end{itemize}

Using this notation, \eqref{genmg} transfers into
\begin{equation}\label{genmgtilde}
\exists\;A\ge 1\;\forall\;p\in\NN_{>0}:\;\;\;\mu_p\le A(\widetilde{M}^d_p)^{\frac{1}{p}}.
\end{equation}
Moreover, this notation should be compared with \eqref{transform} and \eqref{transformM}; see also \cite[Prop. 5.1, $(I)(i)\Rightarrow(ii)$]{modgrowthstrange}: Given $\omega\in\hyperlink{omset0}{\mathcal{W}_0}$ with associated weight matrix $\mathcal{M}_{\omega}=\{\mathbf{W}^{(\ell)}: \ell>0\}$, then \eqref{transform} precisely gives
\begin{equation}\label{transformfct}
\forall\;\ell>0\;\forall\;a\in\NN_{>0}:\;\;\;\mathbf{W}^{(a\ell)}=\widetilde{\mathbf{W}^{(\ell)}}^a,
\end{equation}
and
\begin{equation}\label{tildeweightsequences2gen}
\forall\;x>0\;\forall\;c,p\in\NN_{>0}:\;\;\;\vartheta^{(cx)}_p:=\frac{W^{(cx)}_p}{W^{(cx)}_{p-1}}=\left(\frac{W^{(x)}_{cp}}{W^{(x)}_{cp-c}}\right)^{1/c}=(\vartheta^{(x)}_{cp-c+1}\cdots\vartheta^{(x)}_{cp})^{1/c}.
\end{equation}
If $\omega=\omega_{\mathbf{M}}$ for $\mathbf{M}\in\hyperlink{LCset}{\mathcal{LC}}$, then
$$\mathbf{M}^{(a)}=\widetilde{\mathbf{M}^{(1)}}^a=\widetilde{\mathbf{M}}^{a}$$
and finally, by \cite[Lemma 6.5, $(6.7)$]{PTTvsmatrix}, we get the following which should be compared with \eqref{goodequivalenceclassic}:
\begin{equation}\label{PTT67}
\forall\;a\in\NN_{>0}\;\exists\;D\ge 1\;\forall\;t\ge 0:\;\;\;\omega_{\widetilde{\mathbf{M}}^a}(t)\le\frac{1}{a}\omega_{\mathbf{M}}(t)\le 2\omega_{\widetilde{\mathbf{M}}^a}(t)+D.
\end{equation}

Using these preparations we show the next two technical observations which are in the spirit of \cite[Prop. 5.1]{modgrowthstrange}. First, we verify how the \emph{moderate growth index} is preserved within a given associated weight matrix when replacing $\mathbf{W}^{(x)}$ by $\mathbf{W}^{(cx)}$, $x>0$ and $c\in\NN_{>0}$ arbitrary.

\begin{lemma}\label{technlemma1}
Let $\omega\in\hyperlink{omset0}{\mathcal{W}_0}$ be given and $\mathcal{M}_{\omega}=\{\mathbf{W}^{(\ell)}: \ell>0\}$ the associated weight matrix. Then
\begin{equation}\label{technlemma1equ}
\forall\;x>0\;\forall\;c\in\NN_{>0}:\;\;\;g(\mathbf{W}^{(cx)})\le g(\mathbf{W}^{(x)})\le 4g(\mathbf{W}^{(cx)}).
\end{equation}
\end{lemma}

\demo{Proof}
First, by \eqref{tildeweightsequences2gen} we have
\begin{equation}\label{technlemma3equ}
\forall\;x>0\;\forall\;c,p\in\NN_{>0}:\;\;\;\vartheta^{(x)}_{c(p-1)}\le\vartheta^{(x)}_{cp-c+1}\le\vartheta^{(cx)}_p\le\vartheta^{(x)}_{cp}.
\end{equation}
If now $g(\mathbf{W}^{(x)})=d$ for some $x>0$, arbitrary but from now on fixed, then by the last estimate in \eqref{technlemma3equ} and by \eqref{transform} we get for any $c,p\in\NN_{>0}$:
$$\vartheta^{(cx)}_p\le\vartheta^{(x)}_{cp}\le A(W^{(x)}_{dcp})^{\frac{1}{dcp}}=A(W^{(cx)}_{dp})^{\frac{1}{dp}}.$$
So \eqref{genmg} is verified for $\mathbf{W}^{(cx)}$ with the same choices for $A$ and $d$ and thus $g(\mathbf{W}^{(cx)})\le d$ for all $c\in\NN_{>0}$.

Assume $g(\mathbf{W}^{(cx)})=d$, then by the first and second estimate in \eqref{technlemma3equ}, since $dcp\le 2dc(p-1)\Leftrightarrow 2dc\le dcp\Leftrightarrow 2\le p$ and since $p\mapsto(W^{(x)}_p)^{1/p}$ is nondecreasing, we estimate as follows for all $c,p\in\NN_{>0}$, $p\ge 2$:
$$\vartheta^{(x)}_{c(p-1)}\le\vartheta^{(cx)}_p\le A(W^{(cx)}_{dp})^{\frac{1}{dp}}=A(W^{(x)}_{dcp})^{\frac{1}{dcp}}\le A(W^{(x)}_{2dc(p-1)})^{\frac{1}{2dc(p-1)}}\le A(W^{(x)}_{4dc(p-1)})^{\frac{1}{4dc(p-1)}}.$$
This verifies the desired estimate for $\mathbf{W}^{(x)}$ with parameters $2d$ and the same $A$ for all $p\in\NN$ satisfying $p=c(q-1)$ for some $q\in\NN_{>0}$, $q\ge 2$. Now, let $p\in\NN$ be such that $c(q-1)<p<cq$ for some $q\in\NN_{>0}$, $q\ge 2$. Then $2dcq\le 4dc(q-1)\Leftrightarrow 4dc\le 2dcq\Leftrightarrow 2\le q$, hence by the above estimate applied to $q$ instead of $p-1$:
$$\vartheta^{(x)}_p\le\vartheta^{(x)}_{cq}\le A(W^{(x)}_{2dcq})^{\frac{1}{2dcq}}\le A(W^{(x)}_{4dc(q-1)})^{\frac{1}{4dc(q-1)}}\le A(W^{(x)}_{4dp})^{\frac{1}{4dp}}.$$
This verifies the desired estimate for $\mathbf{W}^{(x)}$ with parameters $4d$ and the same $A$ for all such $p$ under consideration.

And finally, if $1<p<c$, then as verified before $\vartheta^{(x)}_p\le\vartheta^{(x)}_c\le A(W^{(x)}_{2dc})^{\frac{1}{2dc}}=:A_1$ and $A_1\ge A$.

Summarizing, \eqref{genmg} for $\mathbf{W}^{(x)}$ is shown for the choices $A_1$ and $4d$ and hence $g(\mathbf{W}^{(x)})\le 4d$ as desired.
\qed\enddemo

\begin{remark}
\emph{Lemma \ref{technlemma1} holds, in particular, if $\omega=\omega_{\mathbf{M}}$ for some $\mathbf{M}\in\hyperlink{LCset}{\mathcal{LC}}$, and in this situation the conclusion follows for $\mathbf{M}=\mathbf{W}^{(1)}$; see \eqref{transformM}.}

\emph{Next note that \eqref{technlemma1equ} is equivalent to}
\begin{equation}\label{technlemma1equvar}
\forall\;y>0\;\forall\;d\in\NN_{>0}:\;\;\;g(\mathbf{W}^{(y)})\le g(\mathbf{W}^{(y/d)})\le 4g(\mathbf{W}^{(y)}).
\end{equation}
\emph{Indeed, this equivalence holds by the correspondences $c=d$ and $y=xc$. On the other hand, \eqref{technlemma1equvar} can be shown directly analogously as \eqref{technlemma1equ} by involving \cite[$(4.10)$, $(5.2)$]{modgrowthstrange} and hence using the estimates}
\begin{equation*}\label{technlemma3equ1}
\forall\;x>0\;\forall\;c,p\in\NN_{>0}:\;\;\;\vartheta^{(x/c)}_{c(p-1)}\le\vartheta^{(x)}_p=(\vartheta^{(x/c)}_{cp}\cdots\vartheta^{(x/c)}_{cp-c+1})^{1/c}\le\vartheta^{(x/c)}_{cp}.
\end{equation*}
\emph{Summarizing, if for some index $x_0>0$ one has that $\mathbf{W}^{(x_0)}$ satisfies \eqref{genmg} with some $d\in\NN_{>0}$, then \eqref{genmg} also holds for any $\mathbf{W}^{(cx_0)}$ and $\mathbf{W}^{(x_0/c)}$ with a precise modification of the parameter $d$.}
\end{remark}

Analogously one can show the following useful Lemma; we verify that the growth comparison between sequences $\mathbf{M}$ and $\mathbf{N}$ can be transferred to their associated weight matrices $\mathcal{M}_{\omega_{\mathbf{M}}}$ and $\mathcal{M}_{\omega_{\mathbf{N}}}$ in a precise way.

\begin{lemma}\label{sequencematrixequivlemma}
Let $\mathbf{M},\mathbf{N}\in\hyperlink{LCset}{\mathcal{LC}}$ be given and consider the matrices $\mathcal{M}_{\omega_{\mathbf{M}}}=\{\mathbf{M}^{(x)}: x>0\}$ and $\mathcal{M}_{\omega_{\mathbf{N}}}=\{\mathbf{N}^{(x)}: x>0\}$. If $\mathbf{M}\hyperlink{preceq}{\preceq}\mathbf{N}$, then
$$\forall\;x\in\NN_{>0}:\;\;\;\mathbf{M}^{(x)}\hyperlink{preceq}{\preceq}\mathbf{N}^{(x)},\hspace{15pt}\mathbf{M}^{(\frac{1}{2x})}\hyperlink{preceq}{\preceq}\mathbf{N}^{(\frac{1}{x})},$$
consequently $\mathcal{M}_{\omega_{\mathbf{M}}}\hyperlink{Mroumpreceq}{\{\preceq\}}\mathcal{M}_{\omega_{\mathbf{N}}}$ and $\mathcal{M}_{\omega_{\mathbf{M}}}\hyperlink{Mbeurpreceq}{(\preceq)}\mathcal{M}_{\omega_{\mathbf{N}}}$ is valid. More precisely, it holds that
\begin{equation}\label{sequencematrixequivlemmaequ}
\exists\;C\ge 1\;\forall\;x\in\NN_{>0}\;\forall\;p\in\NN:\;\;\;M^{(x)}_p\le C^pN^{(x)}_p,
\end{equation}
and
\begin{equation}\label{sequencematrixequivlemmaequ1}
\forall\;x\in\NN_{>0}\;\exists\;C_x\ge 1\;\forall\;p\in\NN:\;\;\;M^{(\frac{1}{2x})}_p\le C_x^pN^{(\frac{1}{x})}_p.
\end{equation}
Here, $C$ in \eqref{sequencematrixequivlemmaequ} can be chosen to be a constant such that $M_p\le C^pN_p$ holds for all $p\in\NN$ and  in \eqref{sequencematrixequivlemmaequ1} set $C_x:=CM^{(\frac{1}{x})}_{2x}$.
\end{lemma}

In particular, equivalent sequences belonging to the set \hyperlink{LCset}{$\mathcal{LC}$} generate both R- and B-equivalent associated weight matrices.

\demo{Proof}
Since $M_0=1=N_0$ relation $\mathbf{M}\hyperlink{preceq}{\preceq}\mathbf{N}$ precisely means $M_p\le C^pN_p$ for some $C\ge 1$ and all $p\in\NN$. Thus \eqref{transform} and \eqref{transformM} give
$$\forall\;x\in\NN_{>0}\;\forall\;p\in\NN:\;\;\;M^{(x)}_p=(M^{(1)}_{xp})^{\frac{1}{x}}=(M_{xp})^{\frac{1}{x}}\le(C^{xp}N_{xp})^{\frac{1}{x}}=C^pN^{(x)}_p,$$
hence \eqref{sequencematrixequivlemmaequ} is verified. Concerning \eqref{sequencematrixequivlemmaequ1}, first by \eqref{transform1} and \eqref{transformM}
$$\forall\;x\in\NN_{>0}\;\forall\;p\in\NN:\;\;\;M^{(\frac{1}{x})}_{xp}=(M^{(1)}_p)^x=(M_p)^x\le C^{xp}(N_p)^x=C^{xp}N^{(\frac{1}{x})}_{xp}.$$
Thus, for any fixed $x\in\NN_{>0}$, the growth relation for having $\mathbf{M}^{(\frac{1}{x})}\hyperlink{preceq}{\preceq}\mathbf{N}^{(\frac{1}{x})}$ holds with constant $C_x:=C$ and for all values $xp$, $p\in\NN$. Then let $q\in\NN$ be such that $2xp<q<2x(p+1)$ for some $p\in\NN$. Since each sequence under consideration is nondecreasing and by involving \eqref{newmoderategrowth} we estimate as follows:
\begin{align*}
M^{(\frac{1}{2x})}_{q}&\le M^{(\frac{1}{2x})}_{2x(p+1)}\le M^{(\frac{1}{x})}_{2xp}M^{(\frac{1}{x})}_{2x}=(M^{(1)}_{2p})^xM^{(\frac{1}{x})}_{2x}=(M_{2p})^xM^{(\frac{1}{x})}_{2x}
\\&
\le C^{2px}M^{(\frac{1}{x})}_{2x}(N_{2p})^x=C^{2px}M^{(\frac{1}{x})}_{2x}N^{(\frac{1}{x})}_{2xp}\le C^{q}M^{(\frac{1}{x})}_{2x}N^{(\frac{1}{x})}_{q}.
\end{align*}
Moreover, $M^{(\frac{1}{2x})}_p\le M_p^{(\frac{1}{x})}$ for any $p\in\NN$ by the order of the sequences and so we have verified \eqref{sequencematrixequivlemmaequ1}.
\qed\enddemo

\emph{Note:} As mentioned in \cite[Rem. 5.2 $(ii)$]{modgrowthstrange}, in order to define the corresponding weighted spaces, by R- resp. B-equivalence and by the fact that weighted spaces are naturally invariant when replacing a weight (matrix) by some equivalent one, for the Roumieu-type it is sufficient to consider the (sub-)matrix $\{\mathbf{M}^{(cx)}: c\in\NN_{>0}\}$ and for the Beurling-type $\{\mathbf{M}^{(cx)}: c^{-1}\in\NN_{>0}\}$ . Here, the index $x>0$ is arbitrary but fixed. Therefore, Lemmas \ref{technlemma1} and \ref{sequencematrixequivlemma} provide sufficient information for the whole weighted structure.

\begin{remark}
\emph{Concerning Lemma \ref{sequencematrixequivlemma} we comment:}
\begin{itemize}
\item[$(*)$] \emph{In general, the converse implication in this result is unclear: $\mathcal{M}_{\omega_{\mathbf{M}}}\hyperlink{Mroumpreceq}{\{\preceq\}}\mathcal{M}_{\omega_{\mathbf{N}}}$ implies
	$$\exists\;C\ge 1\;\exists\;x\in\NN_{>0}\;\forall\;p\in\NN:\;\;\;M_p=M^{(1)}_p\le C^pN^{(x)}_p=C^p(N_{xp})^{\frac{1}{x}},$$
	and so $M_p^x\le C^{xp}N_{xp}$. Similarly, $\mathcal{M}_{\omega_{\mathbf{M}}}\hyperlink{Mbeurpreceq}{(\preceq)}\mathcal{M}_{\omega_{\mathbf{N}}}$ implies
	$$\exists\;C\ge 1\;\exists\;x\in\NN_{>0}\;\forall\;p\in\NN:\;\;\;M_p^x=(M^{(1)}_p)^x=M^{(\frac{1}{x})}_{xp}\le C^{xp}N^{(1)}_{xp}=C^{xp}N_{xp}.$$
	Thus, in order to verify $\mathbf{M}\hyperlink{preceq}{\preceq}\mathbf{N}$, either $\mathbf{M}$ or $\mathbf{N}$ is required to satisfy moderate growth; see the proof of \cite[Thm. 4.4, $(iii)\Rightarrow(ii)$]{orliczpaper}.}

\item[$(*)$] \emph{Lemma \ref{sequencematrixequivlemma} strengthens one part of the conclusions in \cite[Lemma 4.9 $(i)$]{modgrowthstrange}: There it has been shown that for equivalent sequences $\mathbf{M},\mathbf{N}\in\hyperlink{LCset}{\mathcal{LC}}$ under the assumption that (at least one) sequence is \emph{strongly nonquasianalytic,} see (M.3) in \cite{Komatsu73} and e.g. $(\beta_1)$, $(\gamma_1)$ on \cite[p. 5]{modgrowthstrange}, we have both R- and B-equivalence between the associated weight matrices. Lemma \ref{sequencematrixequivlemma} shows that this assumption and the arguments involving \eqref{assostrongnq} are superfluous.}

\emph{However, in \cite[Lemma 4.9 $(i)$]{modgrowthstrange} even $\omega_{\mathbf{M}}\hyperlink{sim}{\sim}\omega_{\mathbf{N}}$ was shown and which implies stronger estimates than \eqref{sequencematrixequivlemmaequ} and \eqref{sequencematrixequivlemmaequ1}: the exponential growth factors (in the variable $p$) can be replaced by a constant; see also e.g. \eqref{mainmgthm3equ2}.}
\end{itemize}
\end{remark}

Next, given $\mathbf{M},\mathbf{N}\in\RR_{>0}^{\NN}$, then define the \emph{convolved sequence} $\mathbf{M}\star\mathbf{N}$ by
\begin{equation}\label{convolvesequ}
M\star N_p:=\min_{0\le q\le p}M_qN_{p-q},\;\,\;p\in\NN;
\end{equation}
see \cite[$(3.15)$]{Komatsu73} and \cite[$(4.1)$]{weightedentireinclusion1}. This notation has been crucially used in \cite[Sect. 4]{weightedentireinclusion1} and we summarize some properties needed in this work; see also \cite[Sect. 4.2, Rem. 4.3 \& 4.4]{weightedentireinclusion1}:

\begin{itemize}
\item[$(a)$] Obviously $\mathbf{M}\star\mathbf{N}=\mathbf{N}\star\mathbf{M}$ and $M\star N_0=1$ provided that $M_0=N_0=1$.

\item[$(b)$] For all $p\in\NN$ we have $M\star N_p\le\min\{M_0N_p,N_0M_p\}$. So, if in addition $M_0=N_0=1$, then $\mathbf{M}\star\mathbf{N}\le\min\{\mathbf{M},\mathbf{N}\}$ is valid.

\item[$(c)$] For any $\mathbf{M}\in\RR_{>0}^{\NN}$ we have that $\mathbf{M}$ and $\mathbf{M}\star\mathbf{M}$ are equivalent if and only if $\mathbf{M}$ satisfies \hyperlink{mg}{$(\on{mg})$}.

\item[$(d)$] In \cite[Lemma 3.5]{Komatsu73} for log-convex $\mathbf{M},\mathbf{N}$ (resp. $\mathbf{M},\mathbf{N}\in\hyperlink{LCset}{\mathcal{LC}}$) it has been shown that:

\begin{itemize}
\item[$(*)$] $\mathbf{M}\star\mathbf{N}$ is also log-convex (resp. $\mathbf{M}\star\mathbf{N}\in\hyperlink{LCset}{\mathcal{LC}}$).

\item[$(*)$] The corresponding quotient sequence $\mu\star\nu$ is obtained when rearranging resp. ordering the sequences $\mu$ and $\nu$ in the order of growth.

\item[$(*)$] This fact yields by \eqref{counting}
$$\forall\;t\ge 0:\;\;\;\Sigma_{\mathbf{M}\star\mathbf{N}}(t)=\Sigma_{\mathbf{M}}(t)+\Sigma_{\mathbf{N}}(t),$$
and so by \eqref{intrepr}
$$\forall\;t\ge 0:\;\;\;\omega_{\mathbf{M}\star\mathbf{N}}(t)=\omega_{\mathbf{M}}(t)+\omega_{\mathbf{N}}(t).$$
\end{itemize}
\end{itemize}

\section{Preparatory results}\label{preparatorysection}
We prove some technical results which are needed for the main statements in the next section. However, the results in this section are formulated in a general (mixed) setting and might have independent importance for further applications.

\subsection{Weight sequence case}
First, let us treat the weight sequence case. Note that not for all results in this section we have to deal necessarily with sequences belonging to the set \hyperlink{LCset}{$\mathcal{LC}$}. More precisely, the property that $\lim_{p\rightarrow+\infty}(M_p)^{1/p}=+\infty$ holds is only required when involving associated weight functions in order to guarantee that $\omega_{\mathbf{M}}$ is well-defined; see \cite[Sect. 2.4 \& 2.5]{regseqpaper} for more details.

\begin{proposition}\label{technlemma2}
Let $\mathbf{M},\mathbf{L}\in\hyperlink{LCset}{\mathcal{LC}}$ be given. Then the following are equivalent:
\begin{itemize}
\item[$(i)$] It holds that
\begin{equation}\label{technlemma1equ0}
\exists\;a\in\NN_{>0}\;\exists\;B\ge 1\;\forall\;t\ge 0:\;\;\;\omega_{\mathbf{L}\star\widetilde{\mathbf{M}}^a}(t)=\omega_{\mathbf{L}}(t)+\omega_{\widetilde{\mathbf{M}}^a}(t)\le\omega_{\mathbf{L}}(Bt)+B.
\end{equation}

\item[$(ii)$] It holds that
\begin{equation}\label{technlemma1equ2}
\exists\;a\in\NN_{>0}\;\exists\;C\ge 1\;\forall\;t\ge 0:\;\;\;\omega_{\mathbf{M}}(t)\le a(\omega_{\mathbf{L}}(Ct)-\omega_{\mathbf{L}}(t))+C.
\end{equation}

\item[$(iii)$] It holds that
$$\exists\;a\in\NN_{>0}:\;\;\mathbf{L}\hyperlink{preceq}{\preceq}\mathbf{L}\star\widetilde{\mathbf{M}}^a;$$
i.e. in view of \eqref{tildeweightsequences} and \eqref{convolvesequ},
\begin{equation}\label{technlemma1equ1}
\exists\;a\in\NN_{>0}\;\exists\;H\ge 1\;\forall\;p\in\NN\;\forall\;0\le q\le p:\;\;\;L_p\le H^pL_q\widetilde{M}^a_{p-q}=H^pL_q(M_{a(p-q)})^{\frac{1}{a}}.
\end{equation}
\end{itemize}
Concerning the parameter $a$, the proof shows the following correspondences: In $(i)$ and $(iii)$ we can take the same $a$; if $(i)$ holds with $a$ then $(ii)$ with $2a$ and if $(ii)$ holds with $a$ then also $(i)$.
\end{proposition}

\demo{Proof}
$(i)\Leftrightarrow(ii)$ By the second estimate in \eqref{PTT67} and \eqref{technlemma1equ0} we obtain for all $t\ge 0$:
$$\omega_{\mathbf{L}}(t)+\frac{1}{2a}\omega_{\mathbf{M}}(t)-\frac{D}{2}\le\omega_{\mathbf{L}}(t)+\omega_{\widetilde{\mathbf{M}}^a}(t)\le\omega_{\mathbf{L}}(Bt)+B,$$
thus $\omega_{\mathbf{M}}(t)\le 2a(\omega_{\mathbf{L}}(Bt)-\omega_{\mathbf{L}}(t))+2aB+Da$, which gives \eqref{technlemma1equ2} with $2a$ and $C:=2aB+Da(\ge B)$.

Conversely, if \eqref{technlemma1equ2} holds then $a\omega_{\mathbf{L}}(t)+\omega_{\mathbf{M}}(t)\le a\omega_{\mathbf{L}}(Ct)+C$ for all $t\ge 0$ and so by the first estimate in \eqref{PTT67}
$$a\omega_{\mathbf{L}}(t)+a\omega_{\widetilde{\mathbf{M}}^a}(t)\le a\omega_{\mathbf{L}}(t)+\omega_{\mathbf{M}}(t)\le a\omega_{\mathbf{L}}(Ct)+C;$$
thus \eqref{technlemma1equ0} is verified with $B:=C$ and the same $a\in\NN_{>0}$.\vspace{6pt}

$(i)\Leftrightarrow(iii)$ We follow \cite[Prop. 3.6]{Komatsu73}; see also \cite[Prop. 9.3.2]{dissertation} and \cite[Prop. 3.6]{testfunctioncharacterization}: On the one hand, by $(i)$ and \eqref{Prop32Komatsu} applied to $\mathbf{L}$ and to $\mathbf{L}\star\widetilde{\mathbf{M}}^a$ we obtain for all $p\in\NN$:
$$L_p=\sup_{t\ge 0}\frac{t^p}{\exp(\omega_{\mathbf{L}}(t))}=\sup_{s\ge 0}\frac{(Bs)^p}{\exp(\omega_{\mathbf{L}}(Bs))}\le e^BB^p\sup_{s\ge 0}\frac{s^p}{\exp(\omega_{\mathbf{L}\star\widetilde{\mathbf{M}}^a}(s))}=e^BB^pL\star\widetilde{M}^a_p;$$
i.e. $\mathbf{L}\hyperlink{preceq}{\preceq}\mathbf{L}\star\widetilde{\mathbf{M}}^a$. Conversely, we can find $a\in\NN_{>0}$ and $C\ge 1$ such that $L_p\le C^pL\star\widetilde{M}^a_p$ for all $p\in\NN$ and hence $\frac{t^p}{L\star\widetilde{M}^a_p}\le\frac{(Ct)^p}{L_p}$ for all $t\ge 0$ and $p\in\NN$. Thus, by definition $\omega_{\mathbf{L}\star\widetilde{\mathbf{M}}^a}(t)\le\omega_{\mathbf{L}}(Ct)$ for all $t\ge 0$ and so \eqref{technlemma1equ0} is verified with this $a$ and $B:=C$.
\qed\enddemo

We apply Proposition \ref{technlemma2} to $\mathbf{M}=\mathbf{L}$.

\begin{corollary}\label{technlemma2cor}
Let $\mathbf{M}\in\hyperlink{LCset}{\mathcal{LC}}$ be given. Then the following are equivalent:
\begin{itemize}
\item[$(i)$] $\omega_{\mathbf{M}}$ satisfies \eqref{technlemma1equ2} for $\mathbf{M}=\mathbf{L}$,

\item[$(ii)$] $\omega_{\mathbf{M}}$ satisfies \hyperlink{om6}{$(\omega_6)$},

\item[$(iii)$] $\mathbf{M}$ satisfies \hyperlink{mg}{$(\on{mg})$}.
\end{itemize}
\end{corollary}

\demo{Proof}
$(ii)\Leftrightarrow(iii)$ is mentioned in Lemma \ref{assoweightomega0}. Concerning $(i)\Leftrightarrow(ii)$ note that \eqref{technlemma1equ2} for $\mathbf{M}=\mathbf{L}$ precisely means $\frac{a+1}{a}\omega_{\mathbf{M}}(t)\le\omega_{\mathbf{M}}(Ct)+\frac{C}{a}$. Let $a\in\NN_{>0}$ be arbitrary and note that $1<\frac{a+1}{a}\le 2$. Hence \hyperlink{om6}{$(\omega_6)$} for $\omega_{\mathbf{M}}$ immediately gives \eqref{technlemma1equ2} with $C:=Ha$.

Conversely, we iterate \eqref{technlemma1equ2} $n$-times, $n\in\NN_{>0}$ chosen (minimal) to guarantee $\left(\frac{a+1}{a}\right)^n\ge 2$; see also Remark \ref{omega6rem}.
\qed\enddemo

We continue by showing the second main technical statement for the weight sequence case.

\begin{proposition}\label{technlemma3}
Let $\mathbf{M},\mathbf{L}\in\RR_{>0}^{\NN}$ be given such that $M_0=1=L_0$ and such that $\mathbf{L}$ is log-convex. Recall also the notation $\lambda_p:=\frac{L_p}{L_{p-1}}$, $p\in\NN_{>0}$ and $\lambda_0:=1$. Then the following are equivalent:
\begin{itemize}
\item[$(i)$] It holds that
\begin{equation}\label{technlemma2equ}
\exists\;a\in\NN_{>0}\;\exists\;B\ge 1\;\forall\;p\in\NN:\;\;\;L_{2p}\le B^pL_p\widetilde{M}^a_p.
\end{equation}

\item[$(ii)$] It holds that
\begin{equation}\label{technlemma1equ1weak}
\exists\;a\in\NN_{>0}\;\exists\;H\ge 1\;\forall\;p,q\in\NN,\;0\le q\le p:\;\;\;L_{p+q}\le H^pL_q\widetilde{M}^a_{p}.
\end{equation}

\item[$(iii)$] It holds that
\begin{equation}\label{crucialcompequmix}
\exists\;a\in\NN_{>0}\;\exists\;A\ge 1\;\forall\;p\in\NN_{>0}:\;\;\;\lambda_p\le A(M_{ap})^{\frac{1}{ap}}=A(\widetilde{M}^a_{p})^{\frac{1}{p}};
\end{equation}
i.e. the mixed variant of \eqref{genmg} between $\mathbf{L}$ and $\mathbf{M}$ resp. the mixed variant of \eqref{mgstrange} between $\mathbf{L}$ and $\widetilde{\mathbf{M}}^a$.
\end{itemize}
Concerning the parameter $a\in\NN_{>0}$ the proof shows the following correspondences: In $(i)$ and $(ii)$ one can choose the same $a$; if $(i)$ holds with $a$ then also $(iii)$; if $(iii)$ is valid for $a$ then $(ii)$ with $2a$.
\end{proposition}

\demo{Proof}
The implication $(ii)\Rightarrow(i)$ is trivial by considering $q=p$. For $(i)\Rightarrow(ii)$ we follow the trick in the proof of \cite[Thm. 1, $3)\Rightarrow 2)$, p. 674-675]{matsumoto} and proceed by induction: Let $p\in\NN_{>0}$ be given, arbitrary but from now on fixed, and assume that \eqref{technlemma1equ1weak} is verified for some $1\le q\le p$. Indeed, the case $p=0=q$ is trivial and $q=p$ is precisely assumption \eqref{technlemma2equ} (with $H:=B$). Then let us verify \eqref{technlemma1equ1weak} for the case $q-1$: Log-convexity for $\mathbf{L}$ gives $\lambda_{p+q}^{-1}=\frac{L_{p+q-1}}{L_{p+q}}\le\frac{L_{q-1}}{L_q}=\lambda_q^{-1}$ and so
$$L_{p+q-1}=\frac{L_{p+q-1}}{L_{p+q}}L_{p+q}\le B^{p}L_q\widetilde{M}^a_{p}\frac{L_{p+q-1}}{L_{p+q}}\le B^{p}L_q\widetilde{M}^a_{p}\frac{L_{q-1}}{L_q}=B^{p}L_{q-1}\widetilde{M}^a_{p};$$
thus \eqref{technlemma1equ1weak} holds with $H:=B$.\vspace{6pt}

$(i)\Rightarrow(iii)$ By assumption and log-convexity for $\mathbf{L}$ we obtain:
$$\exists\;B\ge 1\;\forall\;p\in\NN:\;\;\;B^{p}(M_{ap})^{\frac{1}{a}}\ge\frac{L_{2p}}{L_p}=\lambda_{2p}\cdots\lambda_{p+1}\ge\lambda_p^p.$$
Therefore, \eqref{crucialcompequmix} is verified with $A:=B$ and the same $a$.\vspace{6pt}

$(iii)\Rightarrow(i)$ By assumption
$$\exists\;a\in\NN_{>0}\;\exists\;A\ge 1\;\forall\;p\in\NN_{>0}:\;\;\;\lambda_{2p}^p\le A^p(M_{2ap})^{\frac{1}{2a}},$$
and this estimate is valid for $p=0$ too having the equality $1=1$. Log-convexity for $\mathbf{L}$ implies that $\lambda_{2p}^p\ge\lambda_{2p}\cdots\lambda_{p+1}=\frac{L_{2p}}{L_p}$ and consequently $L_{2p}\le A^pL_p(M_{2ap})^{\frac{1}{2a}}=A^pL_p\widetilde{M}^{2a}_p$ holds; i.e. \eqref{technlemma2equ} with $2a$ and $B:=A$.
\qed\enddemo

\begin{corollary}\label{technlemma23cor}
Let $\mathbf{M},\mathbf{L}\in\hyperlink{LCset}{\mathcal{LC}}$ be given.
\begin{itemize}
\item[$(i)$] Any of the equivalent assertions in Proposition \ref{technlemma2} implies the equivalent assertions listed in Proposition \ref{technlemma3} and the proof shows that if \eqref{technlemma1equ1} holds with some $a$, then \eqref{technlemma2equ} as well.

\item[$(ii)$] When $\mathbf{M}=\mathbf{L}$, then the converse of $(i)$ is not valid in general; in particular the \emph{q-Gevrey sequences} $M_p:=q^{p^2}$, $q>1$, provide a (counter-)example.
\end{itemize}
\end{corollary}

\demo{Proof}
$(i)$ Evaluate \eqref{technlemma1equ1} at all even integers $2p$ and $q=p$, $p\in\NN$ arbitrary, and hence obtain \eqref{technlemma2equ} with the same $a$ and $B:=H^2$.\vspace{6pt}

$(ii)$ This follows by combining the information from Propositions \ref{technlemma2}, \ref{technlemma3} and Corollary \ref{technlemma2cor}: Consider $\mathbf{M}\in\hyperlink{LCset}{\mathcal{LC}}$ satisfying \eqref{genmg}, i.e. \eqref{crucialcompequmix}, for some $a\in\NN_{\ge 2}$ but violating this condition for $a=1$ and so \hyperlink{mg}{$(\on{mg})$} fails. The q-Gevrey sequences provide such an explicit example; see also \cite[Rem. 4.5 $(ii)$]{modgrowthstrange} for details.
\qed\enddemo

\begin{remark}\label{weakseprem}
\emph{\eqref{technlemma1equ1} and \eqref{technlemma1equ1weak} should be compared with the ``weak separativity condition'' (W.S.) in \cite[p. 669]{matsumoto} between the sequences $\mathbf{L}$ and $\widetilde{\mathbf{M}}^a$. Indeed, for given $\mathbf{L}\in\RR_{>0}^{\NN}$ this condition reads as follows:}
\begin{equation}\label{weaksep}
\exists\;A>0\;\exists\;\mathbf{M}\in\RR_{>0}^{\NN}\;\forall\;p,q\in\NN:\;\;\;L_{p+q}\le A^pL_pM_q.
\end{equation}
\emph{\eqref{weaksep} for the particular choice $\mathbf{M}=\widetilde{\mathbf{M}}^a$ obviously implies both \eqref{technlemma1equ1} and \eqref{technlemma2equ}. If $\mathbf{L}$ is log-convex, then \eqref{weaksep} gives $\lambda_p^p\le\lambda_{p+1}\cdots\lambda_{2p}=\frac{L_{2p}}{L_p}\le A^pM_p$ and so $\lambda_p\le A(M_p)^{1/p}$ for all $p\in\NN_{>0}$; i.e. the mixed variant of \eqref{mgstrange} resp. \eqref{genmg} between $\mathbf{L}$ and $\mathbf{M}$.}
\end{remark}

A first consequence of the established techniques is the following result which should be compared with \cite[Prop. 3.2, 3.3 \& 4.2]{modgrowthstrange} and, when $\mathbf{M}=\mathbf{L}$, with the characterization mentioned in \cite[Thm. 3.1]{modgrowthstrange}.

\begin{proposition}\label{technlemma3consequ}
Let $\mathbf{M},\mathbf{L}\in\hyperlink{LCset}{\mathcal{LC}}$ be given. Consider the following assertions:
\begin{itemize}
\item[$(i)$] \eqref{crucialcompequmix} holds; i.e.
$$\exists\;a\in\NN_{>0}\;\exists\;A\ge 1\;\forall\;p\in\NN_{>0}:\;\;\;\lambda_p\le A(M_{ap})^{\frac{1}{ap}}=A(\widetilde{M}^a_{p})^{\frac{1}{p}}.$$

\item[$(ii)$] It holds that
$$\exists\;a\in\NN_{>0}\;\exists\;A\ge 1\;\forall\;p\in\NN:\;\;\;\lambda_{2p}\le A\widetilde{\mu}^{2a}_p.$$

\item[$(iii)$] It holds that
$$\exists\;a\in\NN_{>0}\;\exists\;A\ge 1\;\forall\;t\ge 0:\;\;\;2\Sigma_{\widetilde{\mathbf{M}}^{2a}}(t)\le\Sigma_{\mathbf{L}}(At).$$

\item[$(iv)$] It holds that
$$\exists\;a\in\NN_{>0}\;\exists\;A\ge 1\;\forall\;t\ge 0:\;\;\;2\omega_{\widetilde{\mathbf{M}}^{2a}}(t)\le\omega_{\mathbf{L}}(At).$$

\item[$(v)$] It holds that $\mathbf{L}\hyperlink{preceq}{\preceq}\widetilde{\mathbf{M}}^{2a}\star\widetilde{\mathbf{M}}^{2a}$, so
 $$\exists\;a\in\NN_{>0}\;\exists\;A\ge 1\;\forall\;p,q\in\NN:\;\;\;L_{p+q}\le A^{p+q}\widetilde{M}^{2a}_p\widetilde{M}^{2a}_q.$$

\item[$(vi)$] It holds that
$$\exists\;a\in\NN_{>0}\;\exists\;A\ge 1\;\forall\;p\in\NN:\;\;\;L_{2p}\le A^p(\widetilde{M}^{2a}_p)^2.$$
\end{itemize}
Then one has the following implications resp. equivalences:
$$(i)\Rightarrow(ii)\Leftrightarrow(iii)\Rightarrow(iv)\Leftrightarrow(v)\Leftrightarrow(vi).$$
Finally, concerning $(vi)\Rightarrow(i)$ we can show a partial converse: If $\mathbf{L}$ satisfies \eqref{genmg} for some $d\in\NN_{>0}$ and
$$\exists\;b\in\NN_{>0}\;\exists\;A\ge 1\;\forall\;p\in\NN:\;\;\;L_{2p}\le A^p(\widetilde{M}^{b}_p)^2,$$
then \eqref{crucialcompequmix} is valid with $a:=db$.
\end{proposition}

\demo{Proof}
$(i)\Rightarrow(ii)$ By \eqref{crucialcompequmix}, log-convexity and normalization for $\mathbf{M}$ we get for all $p\in\NN_{>0}$
$$\lambda_{2p}\le A(\widetilde{M}^a_{2p})^{\frac{1}{2p}}=A(M_{2ap})^{\frac{1}{2ap}}=A(\widetilde{M}^{2a}_{p})^{\frac{1}{p}}\le A\widetilde{\mu}^{2a}_p.$$
Moreover this estimate is clear for $p=0$ as well.\vspace{6pt}

$(ii)\Leftrightarrow(iii)$ This is valid by the definition of the counting functions in \eqref{counting}: Indeed, assume $\lambda_{2p}\le A\widetilde{\mu}^{2a}_p$ for all $p\in\NN$ and some $A\ge 1$ and let $t\ge 0$ be given. If $t\ge\widetilde{\mu}^{2a}_1$ then $\widetilde{\mu}^{2a}_p\le t<\widetilde{\mu}^{2a}_{p+1}$ for some $p\in\NN_{>0}$ and $\Sigma_{\mathbf{L}}(At)\ge\Sigma_{\mathbf{L}}(A\widetilde{\mu}^{2a}_p)\ge\Sigma_{\mathbf{L}}(\lambda_{2p})\ge 2p=2\Sigma_{\widetilde{\mathbf{M}}^{2a}}(t)$ follows. If $0\le t<\widetilde{\mu}^{2a}_1$, then $\Sigma_{\widetilde{\mathbf{M}}^{2a}}(t)=0$ and the desired estimate in $(iii)$ is clear.

Conversely, first we evaluate the identity in $(iii)$ at all $t:=\widetilde{\mu}^a_p$ such that $\widetilde{\mu}^{2a}_p<\widetilde{\mu}^{2a}_{p+1}$: Then $2p=2\Sigma_{\widetilde{\mathbf{M}}^{2a}}(\widetilde{\mu}^{2a}_p)\le\Sigma_{\mathbf{L}}(A\widetilde{\mu}^{2a}_p)$ which gives by definition $A\widetilde{\mu}^{2a}_p\ge\lambda_{2p}$. If $p\in\NN_{>0}$ is such that $\widetilde{\mu}^{2a}_p=\dots=\widetilde{\mu}^{2a}_{p+j}<\widetilde{\mu}^{2a}_{p+j+1}$ for some $j\in\NN_{>0}$, then the first step yields $A\widetilde{\mu}^{2a}_{p+j}\ge\lambda_{2(p+j)}$. Since $p\mapsto\lambda_p$ is nondecreasing, this estimate implies $A\widetilde{\mu}^{2a}_{p+i}\ge\lambda_{2(p+i)}$ for all $0\le i\le j-1$, too.\vspace{6pt}

$(iii)\Rightarrow(iv)$ By \eqref{intrepr} and assumption we get for all $t\ge 0$:
$$2\omega_{\widetilde{\mathbf{M}}^{2a}}(t)=\int_0^t\frac{2\Sigma_{\widetilde{\mathbf{M}}^{2a}}(s)}{s}ds\le\int_0^t\frac{\Sigma_{\mathbf{L}}(As)}{s}ds=\int_0^{At}\frac{\Sigma_{\mathbf{L}}(u)}{u}du=\omega_{\mathbf{L}}(At).$$
$(iv)\Leftrightarrow(v)$ This follows by the same proof as $(i)\Leftrightarrow(iii)$ in Proposition \ref{technlemma2}.\vspace{6pt}

$(v)\Rightarrow(vi)$ is clear and $(vi)\Rightarrow(v)$ follows by the proof of \cite[Thm. 9.5.1, $(i)\Leftrightarrow(ii)$]{dissertation} applied to $\mathbf{M}^{(l)}=\mathbf{L}$, $\mathbf{M}^{(n)}=\widetilde{\mathbf{M}}^{2a}$. Note that in \cite{dissertation} we have involved different arguments compared with the induction used in the proof of \cite[Thm. 1]{matsumoto}: This technique cannot be transferred to the mixed setting.\vspace{6pt}

Concerning the supplement, we estimate as follows for all $p\in\NN_{>0}$:
$$\lambda_p\le A_1(L_{dp})^{\frac{1}{dp}}\le A_1(L_{2dp})^{\frac{1}{2dp}}\le A_1A(M_{dbp})^{\frac{1}{dbp}}.$$
\qed\enddemo

The second consequence is the following result:

\begin{theorem}\label{technlemma4}
Let $\mathbf{M}$ be a given log-convex sequence such that $M_0=1$ and satisfying \eqref{genmg} for some $d\in\NN_{>0}$. Then any log-convex sequence $\mathbf{N}$ being equivalent to $\mathbf{M}$ (and with $N_0=1$) has \eqref{genmg} too; i.e. \eqref{genmg} is preserved under equivalence of log-convex sequences.

Indeed, if $\mathbf{M}$ satisfies \eqref{genmg} with $d$, then for any equivalent log-convex $\mathbf{N}$ we can choose $2d$.
\end{theorem}

\emph{Note:} If $\mathbf{M}$ satisfies \eqref{genmg} with $d=1$ then (equivalently) \hyperlink{mg}{$(\on{mg})$} holds and hence, since this property is obviously preserved under equivalence, also each equivalent log-convex sequence $\mathbf{N}$ satisfies \eqref{genmg} with $d=1$.

\demo{Proof}
We apply Proposition \ref{technlemma3} to $\mathbf{M}=\mathbf{L}$ and get \eqref{technlemma2equ} for $\mathbf{M}$ with parameter $a=2d$. Then, by equivalence,
\begin{align*}
&\exists\;C,D\ge 1\;\forall\;p\in\NN:
\\&
N_{2p}\le C^{2p}M_{2p}\le C^{2p}D^pM_p\widetilde{M}^{2d}_p=C^{2p}D^pM_p(M_{2dp})^{\frac{1}{2d}}\le C^{2p}D^pC^pC^pN_p(N_{2dp})^{\frac{1}{2d}},
\end{align*}
which verifies \eqref{technlemma2equ} for $\mathbf{N}$ with parameter $a=2d$ and $B:=C^4D$. Proposition \ref{technlemma3} applied to $\mathbf{N}=\mathbf{L}$ yields \eqref{genmg} with $2d$.
\qed\enddemo

\begin{remark}\label{technlemma4rem}
\emph{By combining Theorem \ref{technlemma4} and Lemma \ref{sequencematrixequivlemma} and in view of \cite[Prop. 4.4]{modgrowthstrange} we have the fact that both quotient-root comparison properties are preserved when switching to a (R- and B-equivalent) weight matrix which is generated by an equivalent weight sequence. Therefore, the conclusion in \cite[Lemma 4.9 $(ii)$]{modgrowthstrange} holds even without part $(b)$ in \cite[Thm. 4.8]{modgrowthstrange}. Thus this part and also \cite[Lemma 4.6]{modgrowthstrange} are becoming superfluous.}
\end{remark}

\subsection{The weight function case}\label{preparatorysectionweightfct}
First, when given $\omega\in\hyperlink{omset0}{\mathcal{W}_0}$ with associated matrix $\mathcal{M}_{\omega}=\{\mathbf{W}^{(\ell)}: \ell>0\}$, then Proposition \ref{technlemma3consequ} can be applied to the sequences $\mathbf{W}^{(\ell)}$. Note that in view of \eqref{transformfct} it follows that all sequences under consideration in this Proposition belong to $\mathcal{M}_{\omega}$. Indeed, we have:
\begin{itemize}
\item[$(*)$] By \cite[Prop. 4.7]{modgrowthstrange} it follows that $\mathcal{M}_{\omega}$ satisfies the desired quotient-root comparison properties if and only if assertion $(i)$ in Proposition \ref{technlemma3consequ} holds for $\mathbf{L}=\mathbf{W}^{(1)}=\mathbf{M}$.

\item[$(*)$] Assertion $(v)$ holds for $\mathbf{L}=\mathbf{W}^{(\ell)}=\mathbf{M}$ and $A=1=a$ by taking into account \eqref{newmoderategrowth} and \eqref{transformfct}.
\end{itemize}

Therefore, a full equivalence in Proposition \ref{technlemma3consequ} fails in general which is illustrated by \cite[Thm. 4.8]{modgrowthstrange}.

Let us study \eqref{technlemma1equ2} between sequences in the corresponding associated weight matrix.

\begin{lemma}\label{technlemma2cor1}
Let $\omega\in\hyperlink{omset0}{\mathcal{W}_0}$ be given and let $\mathcal{M}_{\omega}=\{\mathbf{W}^{(\ell)}: \ell>0\}$ be the associated weight matrix. Consider the following assertions:
\begin{itemize}
\item[$(i)$] $\omega$ satisfies \hyperlink{om6}{$(\omega_6)$}.

\item[$(ii)$] It holds that
$$\exists\;\ell_1,\ell,a>0\;\exists\;C\ge 1\;\forall\;t\ge 0:\;\;\;\omega_{\mathbf{W}^{(\ell)}}(t)\le a(\omega_{\mathbf{W}^{(\ell_1)}}(Ct)-\omega_{\mathbf{W}^{(\ell_1)}}(t))+C.$$

\item[$(iii)$] It holds that
$$\forall\;\ell_1>0\;\exists\;\ell,a>0\;\exists\;C\ge 1\;\forall\;t\ge 0:\;\;\;\omega_{\mathbf{W}^{(\ell)}}(t)\le a(\omega_{\mathbf{W}^{(\ell_1)}}(Ct)-\omega_{\mathbf{W}^{(\ell_1)}}(t))+C.$$

\item[$(iv)$] It holds that
$$\forall\;\ell>0\;\exists\;\ell_1,a>0\;\exists\;C\ge 1\;\forall\;t\ge 0:\;\;\;\omega_{\mathbf{W}^{(\ell)}}(t)\le a(\omega_{\mathbf{W}^{(\ell_1)}}(Ct)-\omega_{\mathbf{W}^{(\ell_1)}}(t))+C.$$
\end{itemize}
Then $(i)\Rightarrow(ii)-(iv)$ and $(iii),(iv)\Rightarrow(ii)$; if in assertion $(ii)$ for the parameters the relation $\ell_1>\ell a$ holds, then $(ii)\Rightarrow(i)$.
\end{lemma}

\emph{Note:} The estimates in $(ii)-(iv)$ are \eqref{technlemma1equ2} for the sequences $\mathbf{W}^{(\ell)}$ and $\mathbf{W}^{(\ell_1)}$ for which the indices are subject to different quantifiers. Equivalently, in $(ii)-(iv)$ we can assume $a\ge 1$ by enlarging this parameter if necessary and also w.l.o.g. $a\in\NN_{>0}$. Moreover, since $\omega_{\mathbf{W}^{(\ell)}}\le\omega_{\mathbf{W}^{(\ell_1)}}$ for any $0<\ell_1\le\ell$, we see that w.l.o.g. $\ell\ge\ell_1$ in $(ii)$ and $(iii)$. And since in this particular matrix the sequences are even ordered w.r.t. their quotients, i.e. $\vartheta^{(\ell_1)}\le\vartheta^{(\ell_2)}$ for any $0<\ell_1\le\ell_2$, by \eqref{counting} we get $\Sigma_{\mathbf{W}^{(\ell_2)}}(t)\le\Sigma_{\mathbf{W}^{(\ell_1)}}(t)$ for all $t\ge 0$ and then \eqref{intrepr} yields for any $C\ge 1$ and $t\ge 0$:
$$\omega_{\mathbf{W}^{(\ell_2)}}(Ct)-\omega_{\mathbf{W}^{(\ell_2)}}(t)=\int_t^{Ct}\frac{\Sigma_{\mathbf{W}^{(\ell_2)}}(s)}{s}ds\le\int_t^{Ct}\frac{\Sigma_{\mathbf{W}^{(\ell_1)}}(s)}{s}ds=\omega_{\mathbf{W}^{(\ell_1)}}(Ct)-\omega_{\mathbf{W}^{(\ell_1)}}(t).$$
Therefore, also in $(iv)$ w.l.o.g. we can assume $\ell\ge\ell_1$.

\demo{Proof}
The implications $(iii),(iv)\Rightarrow(ii)$ are trivial and we proceed with $(i)\Rightarrow(ii),(iii),(iv)$. First, by combining \hyperlink{om6}{$(\omega_6)$} and \eqref{goodequivalenceclassic} we get the following estimate (recall Remark \ref{omega6rem}):
\begin{align*}
&\forall\;b>0\;\exists\;H\ge 1\;\forall\;\ell_1>0\;\exists\;D_{\ell_1}\ge 1\;\forall\;\ell>0\;\forall\;t\ge 0:
\\&
\ell\omega_{\mathbf{W}^{(\ell)}}(t)+b\ell_1\omega_{\mathbf{W}^{(\ell_1)}}(t)\le(1+b)\omega(t)\le\omega(Ht)+H\le 2\ell_1\omega_{\mathbf{W}^{(\ell_1)}}(Ht)+D_{\ell_1}+H.
\end{align*}
We apply this to $b:=2$ and so
$$\omega_{\mathbf{W}^{(\ell)}}(t)\le\frac{2\ell_1}{\ell}(\omega_{\mathbf{W}^{(\ell_1)}}(Ht)-\omega_{\mathbf{W}^{(\ell_1)}}(t))+\frac{D_{\ell_1}+H}{\ell},$$
which proves all assertions $(ii),(iii),(iv)$ with $a:=\frac{2\ell_1}{\ell}$. (Indeed, note that here $H$ is not depending on the weight matrix indices and hence not on $a$.)\vspace{6pt}

$(ii)\Rightarrow(i)$ We combine the assumption with \eqref{goodequivalenceclassic} and obtain for all $t\ge 0$:
\begin{align*}
&\frac{1}{2\ell}\omega(t)-\frac{D_{\ell}}{2\ell}+\frac{a}{2\ell_1}\omega(t)-\frac{aD_{\ell_1}}{2\ell_1}\le\omega_{\mathbf{W}^{(\ell)}}(t)+a\omega_{\mathbf{W}^{(\ell_1)}}(t)\le a\omega_{\mathbf{W}^{(\ell_1)}}(Ct)+C\le\frac{a}{\ell_1}\omega(Ct)+C.
\end{align*}
Consequently, with $B_{a,C,\ell,\ell_1}:=\frac{C\ell_1}{a}+\frac{D_{\ell}\ell_1}{2\ell a}+\frac{D_{\ell_1}}{2}$, we get
$$\forall\;t\ge 0:\;\;\;\omega(t)\left(\frac{\ell_1+\ell a}{2\ell a}\right)\le\omega(Ct)+B_{a,C,\ell,\ell_1}.$$
If $\ell a\ge\ell_1$, equivalently $\frac{\ell_1+\ell a}{2\ell a}\le 1$, then this estimate is trivial since $C\ge 1$ and $\omega$ is nondecreasing. However, if $\ell_1>\ell a$, then Remark \ref{omega6rem} implies \hyperlink{om6}{$(\omega_6)$} for $\omega$.
\qed\enddemo

\begin{itemize}
\item[$(*)$] If $(ii)$ in Lemma \ref{technlemma2cor1} holds with $\ell_1>\ell a$ and $a\in\NN_{>0}$, then as seen before \hyperlink{om6}{$(\omega_6)$} is valid for $\omega$. In this situation, by Proposition \ref{technlemma2}, \eqref{transformfct} and $(b)$ in Section \ref{auxsection} for the convolved sequence one infers $$\mathbf{W}^{(\ell_1)}\hyperlink{preceq}{\preceq}\mathbf{W}^{(\ell_1)}\star\widetilde{\mathbf{W}^{(\ell)}}^a=\mathbf{W}^{(\ell_1)}\star\mathbf{W}^{(\ell a)}\le\mathbf{W}^{(\ell a)},$$
which is consistent with $(iii)$ in Section \ref{assomatrixsection}.

\item[$(*)$] On the other hand, if
$$\exists\;\ell_1,\ell>0\;\exists\;a\in\NN_{>0},\;\ell_1>2\ell a:\;\;\;\mathbf{W}^{(\ell_1)}\hyperlink{preceq}{\preceq}\mathbf{W}^{(\ell_1)}\star\widetilde{\mathbf{W}^{(\ell)}}^a,$$
then Proposition \ref{technlemma2} yields
$$\exists\;C\ge 1\;\forall\;t\ge 0:\;\;\;\omega_{\mathbf{W}^{(\ell)}}(t)\le 2a(\omega_{\mathbf{W}^{(\ell_1)}}(Ct)-\omega_{\mathbf{W}^{(\ell_1)}}(t))+C,$$
and the arguments for showing $(ii)\Rightarrow(i)$ in Lemma \ref{technlemma2cor1} (applied to $2a$) give \hyperlink{om6}{$(\omega_6)$} for $\omega$.
\end{itemize}

Indeed, this last implication should be compared with the proofs of \cite[Lemma 5.9 $(5.11)$]{compositionpaper} and \cite[Prop. 5.2.1]{dissertation}: The arguments there, verifying \hyperlink{om6}{$(\omega_6)$} for $\omega$, involve the assumption $\mathbf{W}^{(4\ell)}\hyperlink{preceq}{\preceq}\mathbf{W}^{(\ell)}$ for some $\ell>0$ (one chooses for simplicity $\ell=1$). But in view of Remark \ref{omega6rem} it suffices to assume for this implication $\mathbf{W}^{(\ell_1)}\hyperlink{preceq}{\preceq}\mathbf{W}^{(\ell)}$ for some $\ell_1,\ell>0$ satisfying $\ell_1>2\ell$ and so:

\begin{proposition}
Let $\omega\in\hyperlink{omset0}{\mathcal{W}_0}$ be given and let $\mathcal{M}_{\omega}=\{\mathbf{W}^{(\ell)}: \ell>0\}$ be the associated weight matrix. If
$$\exists\;\ell_1,\ell>0,\;\ell_1>2\ell:\;\;\;\mathbf{W}^{(\ell_1)}\hyperlink{preceq}{\preceq}\mathbf{W}^{(\ell)},$$
then $\omega$ satisfies \hyperlink{om6}{$(\omega_6)$}.
\end{proposition}

\section{Main results}\label{mainsection}
In this section we prove the main statements; the weight function setting is then reduced to the weight sequence case. Again, we treat a general (mixed) setting and then restrict to special cases.

\subsection{The weight sequence case}

\begin{lemma}\label{step1}
Let $\mathbf{M},\mathbf{N}\in\hyperlink{LCset}{\mathcal{LC}}$ be given, then the following are equivalent:
\begin{itemize}
\item[$(i)$] It holds that
\begin{equation}\label{auxiweightfct1}
\exists\;C\ge 0\;\exists\;B\ge 1\;\exists\;\ell>0\;\forall\;t\ge 0:\;\;\;2\ell\omega_{\mathbf{N}}(t)\le\omega_{\mathbf{M}}((Bt)^{\ell})+C.
\end{equation}
\item[$(ii)$] It holds that
\begin{equation}\label{auxiweightfct2}
\exists\;A,A_1\ge 1\;\exists\;\ell>0\;\forall\;p\in\NN:\;\;\;M_{2p}\le A_1A^{2p\ell}N_p^{2\ell}.
\end{equation}
\end{itemize}
The parameter $\ell$ can be chosen to be the same in \eqref{auxiweightfct1} and \eqref{auxiweightfct2}; moreover we have the correspondences $A_1=e^C$ and $B=A$.
\end{lemma}

\demo{Proof}
Indeed, this equivalence holds similarly as the proof of \cite[Thm. 7.7, $(iv)\Leftrightarrow(v)$]{orliczpaper}; the techniques are following the ideas in $(ii)\Leftrightarrow(iii)$ in Proposition \ref{technlemma2}. Note that in \cite{orliczpaper} we have $\ell>1$, $\mathbf{M}=\mathbf{N}$ and the factor $A^{2p\ell}$ has not been considered.\vspace{6pt}

$(i)\Rightarrow(ii)$ By using \eqref{Prop32Komatsu} we get for all $p\in\NN$:
\begin{align*}
M_{2p}&=\sup_{t\ge 0}\frac{t^{2p}}{\exp(\omega_{\mathbf{M}}(t))}=\sup_{s\ge 0}\frac{(Bs)^{2p\ell}}{\exp(\omega_{\mathbf{M}}((Bs)^{\ell}))}\le e^C\sup_{s\ge 0}\frac{(Bs)^{2p\ell}}{\exp(2\ell\omega_{\mathbf{N}}(t))}
\\&
=e^CB^{2p\ell}\left(\sup_{s\ge 0}\frac{s^p}{\exp(\omega_{\mathbf{N}}(s))}\right)^{2\ell}=e^CB^{2p\ell}N_p^{2\ell},
\end{align*}
so \eqref{auxiweightfct2} is verified with $A_1:=e^C$, $A:=B$ and the same $\ell>0$.\vspace{6pt}

$(ii)\Rightarrow(i)$ We have $\frac{t^p}{N_p^{\ell}}\le\sqrt{A_1}A^{p\ell}\frac{t^p}{(M_{2p})^{1/2}}$ for all $p\in\NN$ and $t\ge 0$. This yields by definition of the associated weight functions
\begin{equation}\label{step1equ}
\forall\;t\ge 0:\;\;\;\omega_{\mathbf{N}^{\ell}}(t)\le\omega_{\widetilde{\mathbf{M}}^2}(A^{\ell}t)+\frac{\log(A_1)}{2}.
\end{equation}
Take $a=2$ in \eqref{PTT67} and then the first estimate there together with \eqref{powersub} and \eqref{step1equ} imply
$$\forall\;t\ge 0:\;\;\;\ell\omega_{\mathbf{N}}(t^{1/\ell})=\omega_{\mathbf{N}^{\ell}}(t)\le\omega_{\widetilde{\mathbf{M}}^2}(A^{\ell}t)+\frac{\log(A_1)}{2}\le \frac{1}{2}\omega_{\mathbf{M}}(A^{\ell}t)+\frac{\log(A_1)}{2},$$
hence \eqref{auxiweightfct1} is verified with $C:=\log(A_1)$, $B:=A$ and the same $\ell$.
\qed\enddemo

Next we recall \cite[Thm. 5.1]{genLegendreconj}; see also \cite[Lemma 4]{Boman98}:

\begin{theorem}\label{multthm}
Let $\mathbf{M},\mathbf{N}\in\hyperlink{LCset}{\mathcal{LC}}$ be given, then
\begin{equation}\label{multthmequ}
\forall\;t\ge 0:\;\;\;\omega_{\mathbf{M}\mathbf{N}}(t)=\omega_{\mathbf{M}}\check{\star}\omega_{\mathbf{N}}(t),
\end{equation}
with
\begin{equation}\label{uppertransformgen}
\sigma\check{\star}\tau(t):=\inf_{s>0}\{\sigma(s)+\tau(t/s)\},\;\;\;t\in[0,+\infty).
\end{equation}
\end{theorem}

Using this statement we can show the following characterization which is the central result of this work.

\begin{theorem}\label{upptrafothm}
Let $\mathbf{M},\mathbf{N},\mathbf{L}\in\hyperlink{LCset}{\mathcal{LC}}$ be given, then the following are equivalent:
\begin{itemize}
\item[$(i)$] It holds that
\begin{equation}\label{technlemma2equsuper}
\exists\;a\in\NN_{>0}\;\exists\;B\ge 1\;\forall\;p\in\NN:\;\;\;L_{2p}\le B^pN_p\widetilde{M}^a_p.
\end{equation}
\item[$(ii)$] It holds that
\begin{equation}\label{technlemma2equsuperfct}
\exists\;a\in\NN_{>0}\;\exists\;A\ge 1\;\exists\;C\ge 0\;\forall\;t\ge 0:\;\;\;\omega_{\mathbf{N}\widetilde{\mathbf{M}}^a}(t^2)=\omega_{\mathbf{N}}\check{\star}\omega_{\widetilde{\mathbf{M}}^a}(t^2)\le\omega_{\mathbf{L}}(At)+C.
\end{equation}
\end{itemize}
The proof shows that we can take in both estimates the same parameter $a$ and note that when \eqref{technlemma2equsuperfct} holds with some $a$, then also for all $a'\ge a$.
\end{theorem}

\demo{Proof}
We involve Lemma \ref{step1} and the proof illustrates that the parameter $\ell>0$ there does not make any effect.\vspace{6pt}

$(i)\Rightarrow(ii)$ Let $\ell>0$ be arbitrary but from now on fixed, then \eqref{technlemma2equsuper} implies \eqref{auxiweightfct2} for $A_1:=1$, $A:=B^{1/(2\ell)}$, $\mathbf{M}=\mathbf{L}$, and $\mathbf{N}=\mathbf{R}$ with $\mathbf{R}$ denoting the auxiliary sequence given by $R_p:=(N_p\widetilde{M}^a_p)^{1/(2\ell)}$. First, the proof of Lemma \ref{step1} implies
\begin{equation}\label{step2auxequ1}
\forall\;t\ge 0:\;\;\;2\ell\omega_{\mathbf{R}}(t)\le\omega_{\mathbf{L}}(\sqrt{B}t^{\ell}).
\end{equation}
Now observe that $\mathbf{R}\in\hyperlink{LCset}{\mathcal{LC}}$: This is valid because $\mathbf{M},\mathbf{N}\in\hyperlink{LCset}{\mathcal{LC}}$, so $\widetilde{\mathbf{M}}^a\in\hyperlink{LCset}{\mathcal{LC}}$ for any $a\in\NN_{>0}$, see Section \ref{auxsection}. And, moreover, the class of sequences \hyperlink{LCset}{$\mathcal{LC}$} is clearly stable under taking the point-wise product and under taking positive powers of sequences. Then \eqref{powersub} applied to $\frac{1}{2\ell}$ gives
\begin{equation}\label{step2auxequ2}
\forall\;t\ge 0:\;\;\;\omega_{\mathbf{R}}(t)=\omega_{(\mathbf{N}\widetilde{\mathbf{M}}^a)^{1/(2\ell)}}(t)=\frac{1}{2\ell}\omega_{\mathbf{N}\widetilde{\mathbf{M}}^a}(t^{2\ell}).
\end{equation}
In the next step, we crucially apply \eqref{multthmequ} in order to get
\begin{equation}\label{step2auxequ3}
\forall\;t\ge 0:\;\;\;\omega_{\mathbf{N}\widetilde{\mathbf{M}}^a}(t)=\omega_{\mathbf{N}}\check{\star}\omega_{\widetilde{\mathbf{M}}^a}(t).
\end{equation}
Thus, by combining \eqref{step2auxequ1}, \eqref{step2auxequ2} and \eqref{step2auxequ3} with $s:=t^{\ell}$ so far we have shown:
\begin{equation}\label{step2auxequ4}
\forall\;s\ge 0:\;\;\;\omega_{\mathbf{N}}\check{\star}\omega_{\widetilde{\mathbf{M}}^a}(s^2)=2\ell\omega_{\mathbf{R}}(s^{1/\ell})\le\omega_{\mathbf{L}}(\sqrt{B}s).
\end{equation}
Consequently, \eqref{technlemma2equsuperfct} is verified with the same $a$, $A:=\sqrt{B}$ and any $C\ge 0$.\vspace{6pt}

$(ii)\Rightarrow(i)$ Replace $t$ by $t^{\ell}$ in \eqref{technlemma2equsuperfct} and use \eqref{powersub} applied to $\frac{1}{2\ell}$ and \eqref{step2auxequ2} in order to get with $\mathbf{R}:=((N_p\widetilde{M}^a_p)^{1/(2\ell)})_{p\in\NN}$
$$\exists\;A\ge 1\;\exists\;C\ge 0\;\forall\;t\ge 0:\;\;\;2\ell\omega_{\mathbf{R}}(t)=2\ell\omega_{(\mathbf{N}\widetilde{\mathbf{M}}^a)^{1/(2\ell)}}(t)=\omega_{\mathbf{N}\widetilde{\mathbf{M}}^a}(t^{2\ell})\le\omega_{\mathbf{L}}(At^{\ell})+C.$$
Thus we have shown \eqref{auxiweightfct1} with $\mathbf{N}=\mathbf{R}$, $\mathbf{M}=\mathbf{L}$, $B:=A^{1/\ell}$ and the same $C$. Then the proof of Lemma \ref{step1} yields
$$\forall\;p\in\NN:\;\;\;L_{2p}\le e^CA^{2p}(R_p)^{2\ell}=e^CA^{2p}N_p\widetilde{M}^a_p;$$
i.e. \eqref{technlemma2equsuper} is verified with the same $a$ and $B:=e^CA^2$.
\qed\enddemo

Note that \eqref{technlemma2equ} is precisely \eqref{technlemma2equsuper} for the special case $\mathbf{L}=\mathbf{N}$; therefore Theorem \ref{upptrafothm} and Proposition \ref{technlemma3} imply the following result.

\begin{corollary}\label{omegaMcharacterzing}
Let $\mathbf{M},\mathbf{L}\in\hyperlink{LCset}{\mathcal{LC}}$ be given. Then the following are equivalent:
\begin{itemize}
\item[$(i)$] \eqref{technlemma2equ} holds; i.e.
$$\exists\;a\in\NN_{>0}\;\exists\;B\ge 1\;\forall\;p\in\NN:\;\;\;L_{2p}\le B^pL_p\widetilde{M}^a_p=B^pL_p(M_{ap})^{\frac{1}{a}}.$$

\item[$(ii)$] It holds that
$$\exists\;a\in\NN_{>0}\;\exists\;A\ge 1\;\exists\;C\ge 0\;\forall\;t\ge 0:\;\;\;\omega_{\mathbf{L}\widetilde{\mathbf{M}}^a}(t^2)=\omega_{\mathbf{L}}\check{\star}\omega_{\widetilde{\mathbf{M}}^a}(t^2)\le\omega_{\mathbf{L}}(At)+C.$$

\item[$(iii)$] \eqref{crucialcompequmix} holds; i.e.
$$\exists\;a\in\NN_{>0}\;\exists\;A\ge 1\;\forall\;p\in\NN_{>0}:\;\;\;\frac{L_p}{L_{p-1}}=\lambda_p\le A(\widetilde{M}^a_{p})^{\frac{1}{p}}=A(M_{ap})^{\frac{1}{ap}}.$$
\end{itemize}
In $(i)$ and $(ii)$ we can take the same $a$; if $(i)$ holds with $a$ then also $(iii)$; if $(iii)$ is valid for $a$ then $(i)$ with $2a$.
\end{corollary}

Finally, the very special case $\mathbf{M}=\mathbf{N}=\mathbf{L}$ gives the following characterization:

\begin{corollary}\label{omegaMcharacterzing1}
Let $\mathbf{M}\in\hyperlink{LCset}{\mathcal{LC}}$ be given, then the following are equivalent:
\begin{itemize}
\item[$(i)$] $\mathbf{M}$ satisfies \eqref{genmg}.

\item[$(ii)$] The associated weight function satisfies
\begin{equation}\label{omegaMgenmg}
\exists\;a\in\NN_{>0}\;\exists\;A\ge 1\;\exists\;C\ge 0\;\forall\;t\ge 0:\;\;\;\omega_{\mathbf{M}\widetilde{\mathbf{M}}^a}(t^2)=\omega_{\mathbf{M}}\check{\star}\omega_{\widetilde{\mathbf{M}}^a}(t^2)\le\omega_{\mathbf{M}}(At)+C.
\end{equation}
\end{itemize}
If \eqref{genmg} holds with some $d\in\NN_{>0}$, then choose $a:=2d$ in \eqref{omegaMgenmg}; if \eqref{omegaMgenmg} is valid for some $a\in\NN_{>0}$, then choose $d:=a$ in \eqref{genmg}.
\end{corollary}

\begin{remark}\label{omegaMcharacterzing1rem}
\emph{If \eqref{omegaMgenmg} holds with $a=1$, then $\mathbf{M}^a=\mathbf{M}$ and by \eqref{powersub} this estimate gives}
$$\exists\;A\ge 1\;\exists\;C\ge 0\;\forall\;t\ge 0:\;\;\;2\omega_{\mathbf{M}}(t)=\omega_{\mathbf{M}\mathbf{M}}(t^2)\le\omega_{\mathbf{M}}(At)+C;$$
\emph{i.e. \hyperlink{om6}{$(\omega_6)$} for $\omega_{\mathbf{M}}$. And by the known characterization $\mathbf{M}$ satisfies \eqref{genmg} with $d=1$. Indeed, in view of \eqref{powersub} condition \hyperlink{om6}{$(\omega_6)$} is precisely \eqref{omegaMgenmg} with $a=1$.}
\end{remark}

In \eqref{technlemma2equsuperfct} the sequence $\widetilde{\mathbf{M}}^a$ appears and not $\mathbf{M}$ directly but when involving \eqref{PTT67} we can show:

\begin{lemma}\label{upptrafothmcor}
Let $\mathbf{M},\mathbf{N},\mathbf{L}\in\hyperlink{LCset}{\mathcal{LC}}$ be given.
\begin{itemize}
\item[$(i)$] If one has
\begin{equation}\label{technlemma2equsuperfctsuf}
\exists\;A\ge 1\;\exists\;C\ge 0\;\forall\;t\ge 0:\;\;\;\omega_{\mathbf{M}\mathbf{N}}(t^2)=\omega_{\mathbf{N}}\check{\star}\omega_{\mathbf{M}}(t^2)\le\omega_{\mathbf{L}}(At)+C,
\end{equation}
then \eqref{technlemma2equsuperfct} holds for any $a\in\NN_{>0}$ with the same parameters $A,C$.

\item[$(ii)$] If \eqref{technlemma2equsuperfct} holds, then
\begin{equation}\label{technlemma2equsuperfctnec}
\exists\;a\in\NN_{>0}\;\exists\;A\ge 1\;\exists\;B\ge 0\;\forall\;t\ge 0:\;\;\;\omega_{\mathbf{M}\mathbf{N}}(t^2)=\omega_{\mathbf{N}}\check{\star}\omega_{\mathbf{M}}(t^2)\le a\omega_{\mathbf{L}}(At)+B.
\end{equation}
Indeed, if \eqref{technlemma2equsuperfct} is valid with some $a$, then in \eqref{technlemma2equsuperfctnec} it suffices to choose $2a$.
\end{itemize}
\end{lemma}

\demo{Proof}
$(i)$ The first estimate in \eqref{PTT67} and \eqref{uppertransformgen} yield
\begin{align*}
\omega_{\mathbf{N}}\check{\star}\omega_{\widetilde{\mathbf{M}}^a}(t^2)&\le\inf_{s>0}\{\omega_{\mathbf{N}}(s)+\frac{1}{a}\omega_{\mathbf{M}}(t^2/s)\}\le\inf_{s>0}\{\omega_{\mathbf{N}}(s)+\omega_{\mathbf{M}}(t^2/s)\}
\\&
=\omega_{\mathbf{N}}\check{\star}\omega_{\mathbf{M}}(t^2)=\omega_{\mathbf{M}\mathbf{N}}(t^2).
\end{align*}
Hence \eqref{technlemma2equsuperfct} holds with the same parameters as in \eqref{technlemma2equsuperfctsuf}.\vspace{6pt}

$(ii)$ The second estimate in \eqref{PTT67} and \eqref{uppertransformgen} imply for any $t\ge 0$:
\begin{align*}
\omega_{\mathbf{N}}\check{\star}\omega_{\widetilde{\mathbf{M}}^a}(t^2)&\ge\inf_{s>0}\{\omega_{\mathbf{N}}(s)+\frac{1}{2a}\omega_{\mathbf{M}}(t^2/s)\}-\frac{D}{2}\ge\frac{1}{2a}\inf_{s>0}\{\omega_{\mathbf{N}}(s)+\omega_{\mathbf{M}}(t^2/s)\}-\frac{D}{2}
\\&
=\frac{1}{2a}\omega_{\mathbf{N}}\check{\star}\omega_{\mathbf{M}}(t^2)-\frac{D}{2}=\frac{1}{2a}\omega_{\mathbf{M}\mathbf{N}}(t^2)-\frac{D}{2}.
\end{align*}
Combining this with \eqref{technlemma2equsuperfct} gives
$$\exists\;a\in\NN_{>0}\;\exists\;A,D\ge 1\;\exists\;C\ge 0\;\forall\;t\ge 0:\;\;\;\omega_{\mathbf{M}\mathbf{N}}(t^2)\le 2a\omega_{\mathbf{L}}(At)+2aC+aD;$$
i.e. \eqref{technlemma2equsuperfctnec} is verified with $2a$, the same $A$ and $B:=2aC+aD$.
\qed\enddemo

\begin{remark}
\emph{The previous result should be compared with the following comment: When $\mathbf{M}$, $\mathbf{N}$ are log-convex and normalized then
$$\forall\;a\in\NN_{>0}:\;\;\;\mathbf{N}\mathbf{M}\le\mathbf{N}\widetilde{\mathbf{M}}^a\le\widetilde{\mathbf{N}}^a\widetilde{\mathbf{M}}^a,$$
see Section \ref{auxsection}, and consequently the definition of the associated function gives
$$\forall\;t\ge 0:\;\;\;\omega_{\widetilde{\mathbf{M}}^a\widetilde{\mathbf{N}}^a}(t)\le\omega_{\mathbf{N}\widetilde{\mathbf{M}}^a}(t)\le\omega_{\mathbf{N}\mathbf{M}}(t).$$}
\emph{In view of Corollaries \ref{omegaMcharacterzing} and \ref{omegaMcharacterzing1} and $(i)$ in Lemma \ref{upptrafothmcor} a sufficient condition for having \eqref{crucialcompequmix} is}
$$\exists\;A\ge 1\;\exists\;B\ge 0\;\forall\;t\ge 0:\;\;\;\omega_{\mathbf{L}\mathbf{M}}(t^2)=\omega_{\mathbf{L}}\check{\star}\omega_{\mathbf{M}}(t^2)\le\omega_{\mathbf{L}}(At)+B,$$
\emph{and a sufficient condition for having \eqref{genmg} is}
$$\exists\;A\ge 1\;\exists\;B\ge 0\;\forall\;t\ge 0:\;\;\;\omega_{\mathbf{M}^2}(t^2)=\omega_{\mathbf{M}}\check{\star}\omega_{\mathbf{M}}(t^2)\le\omega_{\mathbf{M}}(At)+B.$$
\emph{By taking into account \eqref{powersub} this property turns out to be \hyperlink{om6}{$(\omega_6)$} and hence this is consistent with the known characterization.}
\end{remark}

\subsection{The weight function case}\label{mainsectweightfct}
First, by \cite[Prop. 4.7]{modgrowthstrange} and Corollary \ref{omegaMcharacterzing1} in order to ensure \eqref{rstrange} and/or \eqref{bstrange} it follows that \eqref{omegaMgenmg} has to be satisfied for $\mathbf{M}=\mathbf{W}^{(1)}$. By \eqref{transformfct} it holds that \eqref{omegaMgenmg} is equivalent to
\begin{equation}\label{omegagenmg}
\exists\;a\in\NN_{>0}\;\exists\;A\ge 1\;\exists\;C\ge 0\;\forall\;t\ge 0:\;\;\;\omega_{\mathbf{W}^{(1)}\mathbf{W}^{(a)}}(t^2)=\omega_{\mathbf{W}^{(1)}}\check{\star}\omega_{\mathbf{W}^{(a)}}(t^2)\le\omega_{\mathbf{W}^{(1)}}(At)+C;
\end{equation}
consequently we have:

\begin{corollary}\label{omegaMcharacterzing2}
Let $\omega\in\hyperlink{omset0}{\mathcal{W}_0}$ be given and let $\mathcal{M}_{\omega}=\{\mathbf{W}^{(\ell)}: \ell>0\}$ be the associated weight matrix. Then the following are equivalent:
\begin{itemize}
\item[$(i)$] \eqref{omegagenmg} is valid.

\item[$(ii)$] $\mathcal{M}_{\omega}$ satisfies \eqref{rstrange} and/or \eqref{bstrange}.
\end{itemize}
\end{corollary}

In the following set $\omega_{a}:=\frac{1}{a}\omega$, $a>0$ arbitrary. The second estimate in \eqref{goodequivalenceclassic} implies for all $t\ge 0$ and any $a\ge 1$:
\begin{align*}
\omega_{\mathbf{W}^{(1)}}\check{\star}\omega_{\mathbf{W}^{(a)}}(t^2)&\ge\inf_{s>0}\{\frac{1}{2}\omega(s)+\frac{1}{2a}\omega(t^2/s)\}-\frac{D_1}{2}-\frac{D_{a}}{2a}=\omega_2\check{\star}\omega_{2a}(t^2)-\frac{D_1}{2}-\frac{D_{a}}{2a}
\\&
\ge\frac{1}{2a}\inf_{s>0}\{\omega(s)+\omega(t^2/s)\}-\frac{D_1}{2}-\frac{D_{a}}{2a}=\frac{1}{2a}\omega\check{\star}\omega(t^2)-\frac{D_1}{2}-\frac{D_{a}}{2a}.
\end{align*}
And the first estimate in \eqref{goodequivalenceclassic} implies for all $t\ge 0$ and $a\ge 1$:
\begin{align*}
\omega_{\mathbf{W}^{(1)}}\check{\star}\omega_{\mathbf{W}^{(a)}}(t^2)&\le\inf_{s>0}\{\omega(s)+\frac{1}{a}\omega(t^2/s)\}=\omega\check{\star}\omega_a(t^2)\le\omega\check{\star}\omega(t^2).
\end{align*}
Note that the first estimates in the above computations even hold for all $a>0$. Therefore, \eqref{omegagenmg} and the first estimate in \eqref{goodequivalenceclassic} yield
$$\exists\;a\in\NN_{>0}\;\exists\;A\ge 1\;\exists\;B\ge 0\;\forall\;t\ge 0:\;\;\;\omega\check{\star}\omega(t^2)\le 2a\omega(At)+B$$
with $B:=aD_1+D_{a}+2aC$, and, on the other hand, \eqref{omegagenmg} follows if $\omega$ satisfies
\begin{equation}\label{omegagenmgnew}
\exists\;A\ge 1\;\exists\;C\ge 0\;\forall\;t\ge 0:\;\;\;2\omega\check{\star}\omega(t^2)\le\omega(At)+C.
\end{equation}

\emph{Note:} Indeed \eqref{omegagenmgnew} implies a stronger variant of \eqref{omegagenmg} with quantifiers ``$\exists\;A\ge 1\;\exists\;C\ge 0\;\forall\;a\ge 1\;\forall\;t\ge 0$''. Analogously as in Remark \ref{omegaMcharacterzing1rem}, when \eqref{omegagenmg} holds with $a=1$ then this condition amounts to having \hyperlink{om6}{$(\omega_6)$} for $\omega_{\mathbf{W}^{(1)}}$ and then \cite[Cor. 5.8]{compositionpaper} (see also \cite[Lemma 5.1.3 $(2)$]{dissertation}) implies \hyperlink{om6}{$(\omega_6)$} for $\omega$, too.

Finally, let us prove that \eqref{omegagenmgnew} is already known:

\begin{lemma}\label{om6newcharacterization}
Let $\omega\in\hyperlink{omset0}{\mathcal{W}_0}$ be given. Then \eqref{omegagenmgnew} holds if and only if $\omega$ satisfies \hyperlink{om6}{$(\omega_6)$}.
\end{lemma}
Indeed, in view of Remark \ref{omega6rem} and by inspecting the proof there one infers that \eqref{omegagenmgnew} is also equivalent to requiring $a\omega\check{\star}\omega(t^2)\le\omega(At)+C$ for some/any $a>1$; and $A$ is depending on $a$ when $a$ is subject to an universal quantifier as in $(c)$ in this Remark.

\demo{Proof}
\eqref{omegagenmgnew} implies \hyperlink{om6}{$(\omega_6)$} for $\omega$: We combine \eqref{powersub}, \eqref{multthmequ}, and \eqref{goodequivalenceclassic} and get
\begin{align*}
4\omega_{\mathbf{W}^{(1)}}(t)&=2\omega_{\mathbf{W}^{(1)}\mathbf{W}^{(1)}}(t^2)=2\omega_{\mathbf{W}^{(1)}}\check{\star}\omega_{\mathbf{W}^{(1)}}(t^2)\le 2\omega\check{\star}\omega(t^2)\le\omega(At)+C
\\&
\le 2\omega_{\mathbf{W}^{(1)}}(At)+D_1+C.
\end{align*}
Hence \hyperlink{om6}{$(\omega_6)$} for $\omega_{\mathbf{W}^{(1)}}$ is verified and so \hyperlink{om6}{$(\omega_6)$} for $\omega$ holds as before by \cite[Cor. 5.8]{compositionpaper}.\vspace{6pt}

Conversely, first \hyperlink{om6}{$(\omega_6)$} for $\omega$ implies this property for $\omega_{\mathbf{W}^{(1)}}$ (again by \cite[Cor. 5.8]{compositionpaper}) and \eqref{goodequivalenceclassic} gives
$$\forall\;t\ge 0:\;\;\;\frac{1}{2}\omega\check{\star}\omega(t^2)-D_1\le\omega_{\mathbf{W}^{(1)}}\check{\star}\omega_{\mathbf{W}^{(1)}}(t^2).$$
Then, by involving again \eqref{multthmequ} and \eqref{powersub}, it follows that there exists some $H\ge 1$ such that for all $t\ge 0$:
\begin{align*}
2\omega\check{\star}\omega(t^2)&\le 4\omega_{\mathbf{W}^{(1)}}\check{\star}\omega_{\mathbf{W}^{(1)}}(t^2)+4D_1=4\omega_{\mathbf{W}^{(1)}\mathbf{W}^{(1)}}(t^2)+4D_1
\\&
=8\omega_{\mathbf{W}^{(1)}}(t)+4D_1\le\omega_{\mathbf{W}^{(1)}}(Ht)+H+4D_1\le\omega(Ht)+H+4D_1.
\end{align*}
The second last inequality in this estimate holds by applying Remark \ref{omega6rem} to $\omega_{\mathbf{W}^{(1)}}$ and so \eqref{omegagenmgnew} is verified.
\qed\enddemo

Lemma \ref{om6newcharacterization} and the observations made just before this result are consistent with $(iii)$ in Section \ref{assomatrixsection}: \hyperlink{om6}{$(\omega_6)$} holds if and only if some/each $\mathbf{W}^{(\ell)}$ satisfies \hyperlink{mg}{$(\on{mg})$} and thus \eqref{genmg} with $d=1$.

\begin{remark}
\emph{Let $\omega$ be nondecreasing and satisfying $\lim_{t\rightarrow+\infty}\omega(t)=+\infty$, i.e. $\omega$ is a weight function in the sense of \cite{genLegendreconj} (and of \cite{index}), then we immediately get that \hyperlink{om6}{$(\omega_6)$} for $\omega$ implies:}
$$2\omega\check{\star}\omega(t)=\inf_{s>0}\{2\omega(s)+2\omega(t/s)\}\le\inf_{s>0}\{\omega(Hs)+\omega(Ht/s)\}+2H=\omega\check{\star}\omega(H^2t)+2H.$$
\emph{Note that this implication also holds by \cite[Rem. 2.3, Thm. 3.4 $(ii)$]{genLegendreconj} where the growth index $\overline{\gamma}(\omega)$ from \cite[Sect. 2.4]{genLegendreconj} is involved in the arguments. However, in general the converse implication is not clear and note that for the arguments in Lemma \ref{om6newcharacterization} crucially \eqref{goodequivalenceclassic} is used and for this, more precisely for the second estimate there, the convexity condition \hyperlink{om4}{$(\omega_4)$} is required.}
\end{remark}

For the sake of completeness let us mention that also \hyperlink{om1}{$(\omega_1)$} admits an analogous characterization:

\begin{lemma}\label{om1newcharacterization}
Let $\omega\in\hyperlink{omset0}{\mathcal{W}_0}$ be given. Then $\omega$ satisfies \hyperlink{om1}{$(\omega_1)$} if and only if
$$\exists\;L\ge 1\;\exists\;C\ge 0\;\forall\;t\ge 0:\;\;\;\omega\check{\star}\omega((2t)^2)\le L\omega(t)+C.$$
\end{lemma}

\emph{Note:} Analogously as in Remark \ref{omega6rem} one can require equivalently that $\omega\check{\star}\omega((at)^2)\le L\omega(t)+C$ for some/any $a>1$.

\demo{Proof}
On the one hand, for all $t\ge 0$:
\begin{align*}
\omega_{\mathbf{W}^{(1)}}(2t)&=\frac{1}{2}\omega_{\mathbf{W}^{(1)}\mathbf{W}^{(1)}}((2t)^2)=\frac{1}{2}\omega_{\mathbf{W}^{(1)}}\check{\star}\omega_{\mathbf{W}^{(1)}}((2t)^2)\le \frac{1}{2}\omega\check{\star}\omega((2t)^2)\le\frac{L}{2}\omega(t)+\frac{C}{2}
\\&
\le L\omega_{\mathbf{W}^{(1)}}(t)+\frac{LD_1}{2}+\frac{C}{2},
\end{align*}
which verifies \hyperlink{om1}{$(\omega_1)$} for $\omega_{\mathbf{W}^{(1)}}$ and \hyperlink{om1}{$(\omega_1)$} for $\omega$ follows by \eqref{goodequivalenceclassic} since this property is obviously preserved under equivalence. On the other hand \hyperlink{om1}{$(\omega_1)$} for $\omega$ implies this property for $\omega_{\mathbf{W}^{(1)}}$ (again by \eqref{goodequivalenceclassic}) and so
\begin{align*}
\omega\check{\star}\omega((2t)^2)&\le 2\omega_{\mathbf{W}^{(1)}}\check{\star}\omega_{\mathbf{W}^{(1)}}((2t)^2)+2D_1=2\omega_{\mathbf{W}^{(1)}\mathbf{W}^{(1)}}((2t)^2)+2D_1
\\&
=4\omega_{\mathbf{W}^{(1)}}(2t)+2D_1\le 4L\omega_{\mathbf{W}^{(1)}}(t)+4L+2D_1\le 4L\omega(t)+4L+2D_1.
\end{align*}
\qed\enddemo

\section{Comments on the possible construction of a (counter-)example}\label{counterexsection}
An open question is if the quotient-root comparison property of Roumieu type \eqref{rstrange} and/or of Beurling-type \eqref{bstrange} can be assumed w.l.o.g.; i.e. by switching to an R-equivalent weight matrix satisfying \eqref{rstrange} resp. to an B-equivalent matrix satisfying \eqref{bstrange}. In view of Lemma \ref{sequencematrixequivlemma}, Theorem \ref{technlemma4} and Remark \ref{technlemma4rem} we know that this question has a negative answer within the weight sequence setting; i.e. when we are looking for an R/B-equivalent matrix which is given by the associated weight function expressed by an \emph{equivalent weight sequence.} However, there could exist an equivalent and abstractly given matrix, not being associated with any weight function in the sense of Braun-Meise-Taylor, and satisfying the desired estimates. Abstractly we can prove:

\begin{proposition}\label{mainmgthm3}
Let $\sigma\in\hyperlink{omset0}{\mathcal{W}_0}$ be given with associated weight matrix $\mathcal{M}_{\sigma}=\{\mathbf{S}^{(\ell)}: \ell>0\}$.
\begin{itemize}
\item[$(i)$] Assume that there exists a weight matrix $\mathcal{N}=\{\mathbf{N}^{(x)}: x\in\mathcal{I}\}$ such that $\mathcal{M}_{\sigma}\{\approx\}\mathcal{N}$ and such that

\item[$(ii)$] $\mathcal{N}$ satisfies the quotient-root comparison property of Roumieu type \eqref{rstrange}.
\end{itemize}
Then $\mathbf{S}^{(1)}$ has to satisfy
\begin{equation}\label{mainmgthm3equ}
\exists\;\ell\in\NN_{>0}\;\exists\;B,B_1\ge 1\;\forall\;p\in\NN_{>0}:\;\;\;\frac{S^{(1)}_{p}}{(S^{(1)}_{\ell(p-1)})^{\frac{1}{\ell}}}\le BB_1^p(S^{(1)}_{\ell p})^{\frac{1}{\ell p}}.
\end{equation}
If one assumes instead that $\mathcal{M}_{\sigma}(\approx)\mathcal{N}$ and such that $\mathcal{N}$ satisfies the quotient-root comparison property of Beurling-type \eqref{bstrange}, then $\mathbf{S}^{(1)}$ has to satisfy
\begin{equation}\label{mainmgthm3equ1}
\exists\;\ell\in\NN_{>0}\;\exists\;B,B_1\ge 1\;\forall\;p\in\NN_{>0}:\;\;\;\frac{(S^{(1)}_p)^{\ell}}{S^{(1)}_{\ell p-1}}\le BB_1^{p}(S^{(1)}_{\ell p})^{\frac{1}{\ell p}}.
\end{equation}
\end{proposition}

\emph{Note:} If $\sigma=\omega_{\mathbf{M}}$ for $\mathbf{M}\in\hyperlink{LCset}{\mathcal{LC}}$, then via \eqref{transformM} one has $\mathbf{M}=\mathbf{M}^{(1)}$ and so \eqref{mainmgthm3equ} resp. \eqref{mainmgthm3equ1} has to be satisfied for $\mathbf{M}$.

\demo{Proof}
First, by R-equivalence of the matrices, the point-wise order of the sequences in the matrices and by assumption \eqref{rstrange} for $\mathcal{N}$ we obtain
\begin{align*}
&\forall\;x_1>0\;\exists\;x,y,z>0,\;z\ge y\ge x\ge x_1,\;\exists\;A,C\ge 1\;\forall\;p\in\NN_{>0}:
\\&
\frac{1}{C^{2p}}\frac{S^{(x_1)}_p}{S^{(z)}_{p-1}}\le\frac{N^{(x)}_p}{N^{(x)}_{p-1}}=\nu^{(x)}_p\le A(N^{(y)}_{p})^{\frac{1}{p}}\le AC(S^{(z)}_{p})^{\frac{1}{p}}.
\end{align*}
Then let us apply this to $x_1=1$ and note that by the point-wise order of the sequences in the matrices w.l.o.g. we can assume $z\in\NN_{>0}$. Thus, by \eqref{transform} applied to $\ell=z$ and $x=1$ we have that
\begin{equation}\label{mainmgthm3equ0matrix1}
\exists\;z\in\NN_{>0}\;\exists\;A,C\ge 1\;\forall\;p\in\NN_{>0}:\;\;\;\frac{S^{(1)}_{p}}{(S^{(1)}_{z(p-1)})^{\frac{1}{z}}}\le ACC^{2p}(S^{(1)}_{zp})^{\frac{1}{zp}}.
\end{equation}
This verifies \eqref{mainmgthm3equ} with $\ell:=z$, $B:=AC$ and $B_1:=C^2$.\vspace{6pt}

Similarly, by B-equivalence of the matrices, by assumption \eqref{bstrange} for $\mathcal{N}$, the point-wise order of the sequences in both matrices and recalling the fact that via \eqref{bstrange} we even have that for any $0<x\le y$ it holds that $\nu^{(x)}\le A\nu^{(y)}$ for some $A\ge 1$ (see also \cite[Rem. 4.1 $(ii)$]{modgrowthstrange}), we obtain
\begin{align*}
&\forall\;z>0\;\exists\;y,x,x_1>0,\;x_1\le x\le y\le z,\;\exists\;A,C\ge 1\;\forall\;p\in\NN_{>0}:
\\&
\frac{1}{C^{2p}}\frac{S^{(x_1)}_p}{S^{(z)}_{p-1}}\le\frac{N^{(x)}_p}{N^{(x)}_{p-1}}=\nu^{(x)}_p\le A(N^{(y)}_{p})^{\frac{1}{p}}\le AC(S^{(z)}_{p})^{\frac{1}{p}}.
\end{align*}
We apply this to $z=1$ and by the point-wise order of the sequences in the matrices we can assume w.l.o.g. $x_1=\frac{1}{x_2}$ with some $x_2\in\NN_{>0}$. Then involve \eqref{transform1} applied to $\ell=x_2$ and $x=1$ and consequently the previous estimate applied to all $p=x_2q$, $q\in\NN_{>0}$, yields
$$\exists\;x_2\in\NN_{>0}\;\exists\;A,C\ge 1\;\forall\;q\in\NN_{>0}:\;\;\;\frac{(S^{(1)}_q)^{x_2}}{S^{(1)}_{x_2q-1}}\le ACC^{2x_2q}(S^{(1)}_{x_2q})^{\frac{1}{x_2q}};$$
i.e. \eqref{mainmgthm3equ1} with $\ell:=x_2$, $B:=AC$ and $B_1:=C^{2\ell}$.
\qed\enddemo

Consequently, when applying these techniques, in order to provide a (counter-)example it suffices to construct explicitly $\mathbf{M}\in\hyperlink{LCset}{\mathcal{LC}}$ such that \eqref{mainmgthm3equ} resp. \eqref{mainmgthm3equ1} is violated. However, the following statement shows that such a sequence cannot exist and therefore, in order to answer the conjecture, one has to come up with a different technique.

\begin{lemma}\label{mainmgthm3problemlemma}
For any given $\mathbf{M}\in\hyperlink{LCset}{\mathcal{LC}}$ the estimates \eqref{mainmgthm3equ} and \eqref{mainmgthm3equ1} are valid with $B=B_1=1$ and for any $\ell\in\NN_{>0}$, $\ell\ge 2$.
\end{lemma}

\demo{Proof}
$p\le\ell(p-1)\Leftrightarrow\ell\le p(\ell-1)$ is valid for all $\ell,p\in\NN_{>0}$, $\ell,p\ge 2$. So let $\ell,p\in\NN_{>0}$, $\ell,p\ge 2$ be given and since $p\mapsto(M_p)^{\frac{1}{p}}$ is nondecreasing, first we get $(M_p)^{\frac{1}{p}}\le(M_{\ell(p-1)})^{\frac{1}{\ell(p-1)}}$. And then we have
\begin{align*}
&(M_{\ell(p-1)})^{\frac{1}{\ell(p-1)}}\le(M_{\ell(p-1)})^{\frac{1}{\ell p}}(M_{\ell p})^{\frac{1}{\ell p^2}}\Leftrightarrow(M_{\ell(p-1)})^{\frac{p}{p-1}}\le M_{\ell(p-1)}(M_{\ell p})^{\frac{1}{p}}
\\&
\Leftrightarrow(M_{\ell(p-1)})^{\frac{1}{p-1}}\le(M_{\ell p})^{\frac{1}{p}}\Leftrightarrow(M_{\ell(p-1)})^{\frac{1}{\ell(p-1)}}\le(M_{\ell p})^{\frac{1}{\ell p}},
\end{align*}
and the last inequality is again valid because $p\mapsto(M_p)^{\frac{1}{p}}$ is nondecreasing. Summarizing, so far we have shown
$$\forall\;\ell,p\in\NN_{>0},\;\ell,p\ge 2:\;\;\;M_p\le(M_{\ell(p-1)})^{\frac{1}{\ell}}(M_{\ell p})^{\frac{1}{\ell p}},$$
and, indeed, this estimate is valid for $p=1$ as well: $M_1\le(M_0M_{\ell})^{\frac{1}{\ell}}=(M_{\ell})^{\frac{1}{\ell}}$ is clear for any $\ell\in\NN_{>0}$ because $p\mapsto(M_p)^{\frac{1}{p}}$ is nondecreasing.\vspace{6pt}

Concerning \eqref{mainmgthm3equ1} first note that $M_p\le(M_{\ell p-1})^{\frac{p}{\ell p-1}}$ because $p\le p\ell-1\Leftrightarrow 1\le p(\ell-1)$ which holds for any $p,\ell\in\NN_{>0}$, $\ell\ge 2$. Therefore, $M_p^{\ell}\le(M_{\ell p-1})^{\frac{\ell p}{\ell p-1}}$ and
\begin{align*}
&(M_{\ell p-1})^{\frac{\ell p}{\ell p-1}}\le M_{\ell p-1}(M_{\ell p})^{\frac{1}{\ell p}}\Leftrightarrow(M_{\ell p-1})^{\frac{1}{\ell p-1}}\le(M_{\ell p})^{\frac{1}{\ell p}},
\end{align*}
where the last estimate holds because $0\le\ell p-1\le\ell p$ and $p\mapsto(M_p)^{1/p}$ is nondecreasing. So the conclusion follows.
\qed\enddemo

In view of Lemma \ref{mainmgthm3problemlemma} we close with some more comments on conditions \eqref{mainmgthm3equ}, \eqref{mainmgthm3equ1}:

\begin{itemize}
\item[$(I)$] If $\mathbf{M}\in\hyperlink{LCset}{\mathcal{LC}}$, then since $p\mapsto(M_p)^{1/p}$ is nondecreasing it holds that $M_{p-1}\le(M_{\ell(p-1)})^{\frac{1}{\ell}}$ for any $\ell\in\NN_{>0}$. Therefore, obviously \eqref{genmg} implies \eqref{mainmgthm3equ} with $B:=A$, $B_1:=1$ and $\ell:=d$ and hence this is consistent with $\mathcal{M}_{\sigma}=\mathcal{M}_{\omega_{\mathbf{M}}}$ and \cite[Prop. 4.4]{modgrowthstrange}. Moreover,
    $$\frac{M_p}{M_{p-1}}\ge\frac{M_p^{\ell}}{M_{\ell p-1}}\Leftrightarrow\mu_p\cdots\mu_{\ell p-1}=\frac{M_{\ell p-1}}{M_{p-1}}\ge M_p^{\ell-1},$$
    which is clear for $\ell=1$ and for $\ell\ge 2$ this holds since by log-convexity $\mu_p\cdots\mu_{\ell p-1}\ge\mu_p^{\ell p-p}=\mu^{p(\ell-1)}\ge M_p^{\ell-1}$ because $\mu_p^p\ge\mu_p\cdots\mu_1=M_p$. Thus \eqref{genmg} also implies \eqref{mainmgthm3equ1} with $B:=A$, $B_1:=1$ and $\ell:=d$.

    In particular, if $\mathbf{M}$ has \hyperlink{mg}{$(\on{mg})$}, then we can take in both conditions $\ell=1=B_1$.

\item[$(II)$] By recognizing Lemma \ref{sequencematrixequivlemma}, in \cite[Lemma 4.6]{modgrowthstrange} we have treated a ``special version'' of Proposition \ref{mainmgthm3} which is corresponding to $\sigma=\omega_{\mathbf{N}}$ and $\mathcal{N}=\mathcal{M}_{\omega_{\mathbf{M}}}$. Via \cite[Prop. 4.4]{modgrowthstrange} in this situation $\mathbf{M}$ has to satisfy \eqref{genmg} and any $\mathbf{N}\in\RR_{>0}^{\NN}$ being equivalent to $\mathbf{M}$ has to satisfy \cite[$(4.13)$]{modgrowthstrange} which reads as follows:
    \begin{equation}\label{equ413new}
    \exists\;d\in\NN_{>0}\;\exists\;A,C\ge 1\;\forall\;p\in\NN_{>0}:\;\;\;N_p\le AC^{2p}(N_{dp})^{\frac{1}{dp}}N_{p-1}.
    \end{equation}
    When $\mathbf{N}$ is log-convex and normalized then we get $N_{p-1}\le(N_{\ell(p-1)})^{\frac{1}{\ell}}\Leftrightarrow N_{p-1}^{\ell}\le N_{\ell(p-1)}$ for any $\ell\in\NN_{>0}$ and so condition \eqref{equ413new} implies \eqref{mainmgthm3equ} with $d:=\ell$, $B:=A$ and $B_1:=C^2$. And also \eqref{mainmgthm3equ1} follows with the same choices when repeating the above arguments in $(I)$ for $\mathbf{N}$.

    However, as mentioned before, by Theorem \ref{technlemma4} this case is completely understood and in this sense the study of \eqref{equ413new} is superfluous.

\item[$(III)$] On the other hand one can ask if Proposition \ref{mainmgthm3} simplifies within the weight function setting; i.e. when restricting to $\mathcal{N}=\mathcal{M}_{\omega}$ for some $\omega\in\hyperlink{omset0}{\mathcal{W}_0}$. More precisely, in this situation consider the following:

\begin{itemize}
\item[$(a)$] Let $\sigma\in\hyperlink{omset0}{\mathcal{W}_0}$ and assume that there exists $\omega\in\hyperlink{omset0}{\mathcal{W}_0}$ such that $\sigma\hyperlink{sim}{\sim}\omega$ and such that

\item[$(b)$] $\mathcal{M}_{\omega}=\{\mathbf{W}^{(\ell)}: \ell>0\}$ satisfies the quotient-root comparison property of Roumieu type \eqref{rstrange} and/or of Beurling-type \eqref{bstrange}.
\end{itemize}

First, by the characterization obtained in \cite[Prop. 4.7]{modgrowthstrange} assertion $(b)$ holds if and only if
\begin{equation}\label{mainmgthm3equ0}
\exists\;d\in\NN_{>0}\;\exists\;A\ge 1\;\forall\;p\in\NN_{>0}:\;\;\;\frac{W^{(1)}_p}{W^{(1)}_{p-1}}=\vartheta^{(1)}_p\le A(W^{(1)}_{dp})^{\frac{1}{dp}}.
\end{equation}
Second, we recall (the proof of) \cite[Lemma 5.16]{compositionpaper} even under more general requirements on $\sigma$: Let $\sigma:[0,+\infty)\rightarrow[0,+\infty)$ be nondecreasing such that $\lim_{t\rightarrow+\infty}\sigma(t)=+\infty$, i.e. $\sigma$ is a weight function in the sense of \cite{genLegendreconj} (and of \cite{index}), and that $\sigma\hyperlink{sim}{\sim}\omega$ holds. This relation yields \hyperlink{om3}{$(\omega_3)$} for $\sigma$ and precisely means
$$\exists\;B\ge 1\;\forall\;t\ge 0:\;\;\;\frac{1}{B}\sigma(t)-1\le\omega(t)\le B\sigma(t)+B.$$
By definition we get the same estimates for $\varphi_{\sigma}=\sigma\circ\exp$ and $\varphi_{\omega}=\omega\circ\exp$ even for all $t\in\RR$ and $\RR$ should be considered as the natural set of definition for the Young-conjugate when the weight is not normalized. Thus, for all $s\ge 0$:
\begin{align*}
\varphi^{*}_{\omega}(s)&=\sup_{t\ge 0}\{ts-\varphi_{\omega}(t)\}=\sup_{t\in\RR}\{ts-\varphi_{\omega}(t)\}\ge\sup_{t\in\RR}\{ts-B\varphi_{\sigma}(t)\}-B
\\&
=B\sup_{t\in\RR}\{t(s/B)-\varphi_{\sigma}(t)\}-B=B\varphi^{*}_{\sigma}(s/B)-B,
\end{align*}
and similarly when $\omega$ and $\sigma$ are interchanged. Consequently,
$$\exists\;B\ge 1\;\forall\;s\ge 0:\;\;\;B\varphi^{*}_{\sigma}(s/B)\le\varphi^{*}_{\omega}(s)+B,\;\;\;\;\;B\varphi^{*}_{\omega}(s)\le\varphi^{*}_{\sigma}(Bs)+B,$$
and when evaluating these estimates at $s:=\ell p$, $p\in\NN$ and $\ell>0$ arbitrary, by definition we obtain
\begin{equation}\label{mainmgthm3equ2}
\exists\;B\ge 1\;\forall\;\ell>0\;\forall\;p\in\NN:\;\;\;e^{-B/\ell}S^{(\ell/B)}_p\le W^{(\ell)}_p\le e^{1/\ell}S^{(\ell B)}_p.
\end{equation}
Then apply \eqref{mainmgthm3equ2} to $\ell:=1$ and combining these estimates with \eqref{mainmgthm3equ0} yields
\begin{equation}\label{mainmgthm3equ3}
\exists\;A,B\ge 1\;\exists\;d\in\NN_{>0}\;\forall\;p\in\NN_{>0}:\;\;\;\frac{e^{-B}S^{(1/B)}_p}{eS^{(B)}_{p-1}}\le Ae(S^{(B)}_{dp})^{\frac{1}{dp}}.
\end{equation}
W.l.o.g. one can assume that $B\in\NN_{>0}$. Then by \eqref{transform} and \eqref{transform1}, both applied to $\ell:=B$ and $x:=1$, we have that \eqref{mainmgthm3equ3} implies the following necessary condition (when $p$ is replaced by $Bp$):
\begin{equation}\label{mainmgthm3equ4}
\exists\;A\ge 1\;\exists\;B,d\in\NN_{>0}\;\forall\;p\in\NN_{>0}:\;\;\;\frac{(S^{(1)}_p)^{B}}{(S^{(1)}_{B(Bp-1)})^{\frac{1}{B}}}\le Ae^{B+2}(S^{(1)}_{B^2dp})^{\frac{1}{B^2dp}}.
\end{equation}
This should be compared with \eqref{mainmgthm3equ}, \eqref{mainmgthm3equ1}. But, similarly as before, we verify that \eqref{mainmgthm3equ4} is valid for all $\mathbf{M}\in\hyperlink{LCset}{\mathcal{LC}}$ with any choices for $A\ge 1$, $B\in\NN_{\ge 2}$ and $d\in\NN_{>0}$: First $(M_p)^B\le(M_{B(Bp-1)})^{\frac{p}{Bp-1}}\Leftrightarrow(M_p)^{1/p}\le(M_{B(Bp-1)})^{\frac{1}{B(Bp-1)}}$ holds for all $B\in\NN_{\ge 2}$ and $p\in\NN_{>0}$ since $p\mapsto(M_p)^{1/p}$ is nondecreasing and $p\le 2(2p-1)\Leftrightarrow\frac{2}{3}\le p$. Then
\begin{align*}
&(M_{B(Bp-1)})^{\frac{p}{Bp-1}}\le(M_{B(Bp-1)})^{\frac{1}{B}}(M_{B^2dp})^{\frac{1}{B^2dp}}
\\&
\Leftrightarrow(M_{B(Bp-1)})^{\frac{p}{Bp-1}-\frac{1}{B}}=(M_{B(Bp-1)})^{\frac{1}{B(Bp-1)}}\le(M_{B^2dp})^{\frac{1}{B^2dp}},
\end{align*}
and this last estimate is valid for any $B\in\NN_{>0}$ and $d\in\NN_{>0}$ since $B(Bp-1)\le B^2dp\Leftrightarrow-1\le Bp(d-1)$.

Summarizing, again it does not make sense to try to construct $\sigma\in\hyperlink{omset0}{\mathcal{W}_0}$ such that $\mathbf{S}^{(1)}$ violates \eqref{mainmgthm3equ4}.
\end{itemize}

\bibliographystyle{plain}
\bibliography{Bibliography}

\begin{thebibliography}{10}

\bibitem{Bjorck66}
G.~Bj{\"o}rck.
\newblock Linear partial differential operators and generalized distributions.
\newblock {\em Ark. Mat.}, 6:351--407, 1966.

\bibitem{nuclearglobal2}
C.~Boiti, D.~Jornet, A.~Oliaro, and G.~Schindl.
\newblock Nuclear global spaces of ultradifferentiable functions in the matrix
  weighted setting.
\newblock {\em Banach J. of Math. Anal.}, 15(1):art. no. 14, 2021.

\bibitem{nuclearglobal1}
C.~Boiti, D.~Jornet, A.~Oliaro, and G.~Schindl.
\newblock Nuclearity of rapidly decreasing ultradifferentiable functions and
  time-frequency {A}nalysis.
\newblock {\em Collect. Math.}, 72(1):423--442, 2021.

\bibitem{Boman98}
J.~Boman.
\newblock Uniqueness and non-uniqueness for microanalytic continuation of
  ultradistributions.
\newblock {\em Contemporary Mathematics}, 251:61--82, 2000.

\bibitem{BonetMeiseMelikhov07}
J.~Bonet, R.~Meise, and S.~N. Melikhov.
\newblock A comparison of two different ways to define classes of
  ultradifferentiable functions.
\newblock {\em Bull. Belg. Math. Soc. Simon Stevin}, 14:424--444, 2007.

\bibitem{BraunMeiseTaylor90}
R.~W. Braun, R.~Meise, and B.~A. Taylor.
\newblock Ultradifferentiable functions and {F}ourier analysis.
\newblock {\em Res. Math.}, 17(3-4):206--237, 1990.

\bibitem{CordaroFuerdoes24}
P.~D. Cordaro and S.~F\"{u}rd\"{o}s.
\newblock The {M}etivier inequality and ultradifferentiable hypoellipticity.
\newblock {\em Math. Nachr.}, 297:2517--2531, 2024.

\bibitem{optimalflat23}
J.~Jim\'{e}nez-Garrido, I.~Miguel-Cantero, J.~Sanz, and G.~Schindl.
\newblock Optimal flat functions in {C}arleman-{R}oumieu ultraholomorphic
  classes in sectors.
\newblock {\em Res. Math.}, 78:art. no. 98, 2023.

\bibitem{PTTvsmatrix}
J.~Jim\'{e}nez-Garrido, D.~N. Nenning, and G.~Schindl.
\newblock On generalized definitions of ultradifferentiable classes.
\newblock {\em J. Math. Anal. Appl.}, 526:127260, 2023.

\bibitem{index}
J.~Jim\'{e}nez-Garrido, J.~Sanz, and G.~Schindl.
\newblock Indices of {O}-regular variation for weight functions and weight
  sequences.
\newblock {\em Rev. R. Acad. Cienc. Exactas Fís. Nat. Ser. A Mat. RACSAM},
  113(4):3659--3697, 2019.

\bibitem{sectorialextensions}
J.~Jiménez-Garrido, J.~Sanz, and G.~Schindl.
\newblock Sectorial extensions, via {L}aplace transforms, in ultraholomorphic
  classes defined by weight functions.
\newblock {\em Res. Math.}, 74(27), 2019.

\bibitem{Komatsu73}
H.~Komatsu.
\newblock Ultradistributions. {I}. {S}tructure theorems and a characterization.
\newblock {\em J. Fac. Sci. Univ. Tokyo Sect. IA Math.}, 20:25--105, 1973.

\bibitem{mandelbrojtbook}
S.~Mandelbrojt.
\newblock {\em Séries adhérentes, Régularisation des suites, Applications}.
\newblock Gauthier-Villars, Paris, 1952.
\newblock (in French).

\bibitem{matsumoto}
W.~Matsumoto.
\newblock Characterization of the separativity of ultradifferentiable classes.
\newblock {\em J. Math. Kyoto Univ.}, 24(4):667--678, 1984.

\bibitem{compositionpaper}
A.~Rainer and G.~Schindl.
\newblock Composition in ultradifferentiable classes.
\newblock {\em Studia Math.}, 224(2):97--131, 2014.

\bibitem{whitneyextensionweightmatrix}
A.~Rainer and G.~Schindl.
\newblock Extension of {W}hitney jets of controlled growth.
\newblock {\em Math. Nachr.}, 290(14--15):2356--2374, 2017.

\bibitem{whitneyextensionmixedweightfunctionII}
A.~Rainer and G.~Schindl.
\newblock On the extension of {W}hitney ultrajets, {I}{I}.
\newblock {\em Studia Math.}, 250(3):283--295, 2020.

\bibitem{genLegendreconj}
G.~Schindl.
\newblock Generalized upper and lower {L}egendre conjugates for weight
  functions.
\newblock 2025, available online at \url{https://arxiv.org/pdf/2505.07497.pdf}.

\bibitem{diploma}
G.~Schindl.
\newblock Spaces of smooth functions of {D}enjoy-{C}arleman-type, 2009.
\newblock Diploma Thesis, Universität Wien, available online at
  \url{http://othes.univie.ac.at/7715/1/2009-11-18_0304518.pdf}.

\bibitem{dissertation}
G.~Schindl.
\newblock Exponential laws for classes of {D}enjoy-{C}arleman differentiable
  mappings, 2014.
\newblock PhD Thesis, Universität Wien, available online at
  \url{http://othes.univie.ac.at/32755/1/2014-01-26_0304518.pdf}.

\bibitem{testfunctioncharacterization}
G.~Schindl.
\newblock Characterization of ultradifferentiable test functions defined by
  weight matrices in terms of their {F}ourier transform.
\newblock {\em Note di Matematica}, 36(2):1--35, 2016.

\bibitem{subaddlike}
G.~Schindl.
\newblock On subadditivity-like conditions for associated weight functions.
\newblock {\em Bull. Belg. Math. Soc. Simon Stevin}, 28(3):399--427, 2022.

\bibitem{modgrowthstrange}
G.~Schindl.
\newblock On the equivalence between moderate growth-type conditions in the
  weight matrix setting.
\newblock {\em Note di Matem.}, 42(1):1--35, 2022.

\bibitem{weightedentireinclusion1}
G.~Schindl.
\newblock On inclusion relations between weighted spaces of entire functions.
\newblock {\em Bull. Sci. Math.}, 190:103375, 2024.

\bibitem{orliczpaper}
G.~Schindl.
\newblock On {O}rlicz classes defined in terms of associated weight functions.
\newblock {\em Monatsh. Math.}, 204:919--968, 2024.

\bibitem{regseqpaper}
G.~Schindl.
\newblock On the regularization of sequences and associated weight functions.
\newblock {\em Bull. Belg. Math. Soc. - Simon Stevin}, 31(2):174--210, 2024.

\bibitem{Thilliezdivision}
V.~Thilliez.
\newblock Division by flat ultradifferentiable functions and sectorial
  extensions.
\newblock {\em Res. Math.}, 44:169--188, 2003.

\end{thebibliography}

\end{document}